\pdfoutput=1
\documentclass[final,3p,times]{elsarticle}

\usepackage{amssymb}
\usepackage{mathtools} %
\usepackage{bm}
\usepackage{physics}
\usepackage{dsfont}
\usepackage{subcaption}
\usepackage{multirow}
\usepackage{algorithm}
\usepackage{algpseudocode}
\usepackage{booktabs}
\usepackage[normalem]{ulem}
\usepackage{verbatim}
\usepackage{siunitx}

\newcommand{\verso}[1]{}
\date{March 26, 2025}

\usepackage{xcolor}

\usepackage{pgfplots}
\usepgfplotslibrary{groupplots,dateplot}
\pgfplotsset{compat=newest}
\usepackage{tikz}
\usetikzlibrary{external}

\usetikzlibrary{patterns,shapes.arrows,external,math, matrix}
\definecolor{mplgreen}{RGB}{44,160,44}  %
\definecolor{mplblue}{RGB}{31,119,180}  %
\definecolor{mplred}{RGB}{214,39,40}    %
\definecolor{pvred}{RGB}{180,4,38}      %
\definecolor{cwred}{RGB}{198,51,51}
\definecolor{cwblue}{RGB}{73,97,210}

\definecolor{cblack}{RGB}{0,0,0}

\definecolor{darkgray176}{RGB}{176,176,176}  %
\definecolor{lightgray204}{RGB}{204,204,204} %
\definecolor{darkorange25512714}{RGB}{255,127,14}
\definecolor{forestgreen4416044}{RGB}{44,160,44}
\definecolor{lightgray204}{RGB}{204,204,204}
\definecolor{mediumpurple148103189}{RGB}{148,103,189}
\definecolor{steelblue31119180}{RGB}{31,119,180}

\colorlet{mixed41}{mplblue!90}
\colorlet{mixed42}{mplblue!50}
\colorlet{mixed43}{cwred!50}
\colorlet{mixed44}{cwred!90}
\colorlet{mixed45}{cwred!180}

\colorlet{mixed31}{mplblue!90}
\colorlet{mixed32}{mplblue!50}
\colorlet{mixed33}{cwred!50}

\colorlet{conv21}{mplblue!80}
\colorlet{conv22}{cwred!80}
\colorlet{conv23}{mplgreen!80}

\colorlet{conv31}{mplblue!80}
\colorlet{conv32}{mplgreen!80}
\colorlet{conv33}{cwred!80}
\colorlet{conv34}{darkorange25512714!80}
\colorlet{conv35}{mediumpurple148103189!80}
\tikzset{order guide one/.style={very thick, black, dash pattern=on 5pt off 2pt}}
\tikzset{order guide two/.style={very thick, black, dash pattern=on 5pt off 2pt on \the\pgflinewidth off 2pt}}
\tikzset{order guide three/.style={very thick, black, dash pattern=on 5pt off 2pt on \the\pgflinewidth off 2pt on \the\pgflinewidth off 2pt}}
\tikzset{stationary/.style={very thick, dash pattern=on 5pt off 2pt}}
\tikzset{convergence circle/.style={very thick, mark=*, mark size=3.0, mark options={solid,fill=white,fill opacity=1}}}
\tikzset{convergence square/.style={very thick, mark=square*, mark size=2.5, mark options={solid,fill=white,fill opacity=1}}}
\tikzset{convergence triangle/.style={very thick, mark=triangle*, mark size=4, mark options={solid,fill=white,fill opacity=1}}}
\tikzset{convergence diamond/.style={very thick, mark=diamond*, mark size=4, mark options={solid,fill=white,fill opacity=1}}}
\tikzset{convergence pentagon/.style={very thick, mark=pentagon*, mark size=4, mark options={solid,fill=white,fill opacity=1}}}
\tikzset{stationary/.style={very thick, dash pattern=on 5pt off 2pt}}
\tikzset{time convergence/.style={very thick}}
\tikzset{spectrum1/.style={draw=conv21, fill=conv21, mark=*, only marks, mark size=1.3pt}}
\tikzset{spectrum2/.style={draw=conv22, fill=conv22, mark=*, only marks, mark size=1.3pt}}
\tikzset{spectrum3/.style={draw=conv23, fill=conv23, mark=*, only marks, mark size=1.3pt}}

\pgfplotsset{mixed aspect/.style={height=190pt, width=240pt}}

\usepackage[export]{adjustbox}

\usepackage{hyperref}
\hypersetup{
    colorlinks=true, 
    pdftitle={Gabbard, Paris, van Rees. JCP 2025},
    pdfauthor={James Gabbard, Andrea Paris, Wim M. van Rees}
}

\newcommand{\dbstext}[1]{\ \ \text{#1}\ \ }

\newcommand{\inlineorder}[1]{\mathcal{O}(#1)}

\newcommand{\xv}{\vb{x}} %
\newcommand{\sv}{\vb{s}} %
\newcommand{\kv}{\vb{k}} %

\newcommand{\ReviewerOneStrike}[1]{}
\newcommand{\ReviewerTwoStrike}[1]{}

\journal{Journal of Computational Physics}

\begin{document}
\begin{frontmatter}

\title{A high order multigrid-preconditioned immersed interface solver for the Poisson equation with boundary and interface conditions}

\author[MIT]{James Gabbard\fnref{fn1}}
\author[MIT,ETH]{Andrea Paris\fnref{fn1}}
\author[MIT]{Wim M. van Rees\corref{cor1}}
\cortext[cor1]{Corresponding author. E-mail address: wvanrees@mit.edu  
}
\affiliation[MIT]{
    organization={Department of Mechanical  Engineering, Massachusetts Institute of Technology},%
    addressline={77 Masschusetts Ave.}, 
    city={Cambridge},
    postcode={02139}, 
    state={MA},
    country={United States}
}
\affiliation[ETH]{
    organization={Department of Mechanical and Process Engineering, ETH Zürich},%
    addressline={Rämistrasse 101}, 
    city={8092 Zürich},
    country={Switzerland}
}
\fntext[fn1]{Contributed equally}

\begin{abstract}
This work presents a multigrid preconditioned high order immersed finite difference solver to accurately and efficiently solve the Poisson equation on complex 2D and 3D domains. The solver employs a low order Shortley-Weller multigrid method to precondition a high order matrix-free Krylov subspace solver. The matrix-free approach enables full compatibility with high order IIM discretizations of boundary and interface conditions, as well as high order wavelet-adapted multiresolution grids. Through verification and analysis on 2D domains, we demonstrate the ability of the algorithm to provide high order accurate results to Laplace and Poisson problems with Dirichlet, Neumann, and/or interface jump boundary conditions, all effectively preconditioned using the multigrid method. 
We further show that the proposed method is able to efficiently solve high order discretizations of Laplace and Poisson problems on complex 3D domains using thousands of compute cores and on multiresolution grids. To our knowledge, this work presents the largest problem sizes tackled with high order immersed methods applied to elliptic partial differential equations, and the first high order results on 3D multiresolution adaptive grids. Together, this work paves the way for employing high order immersed methods to a variety of 3D partial differential equations with boundary or interface conditions, including linear and non-linear elasticity problems, the incompressible Navier-Stokes equations, and fluid-structure interactions. 
\end{abstract}


\end{frontmatter}

\section{Introduction}\label{sec:introduction}

Elliptic problems arise in partial differential equations (PDEs) across disciplines, notably in incompressible continuum mechanics as well as various equilibrium problems in heat transfer, electrostatics, and magnetostatics. To solve such problems within complex domains, immersed discretizations on regular, structured grids \cite{Peskin1972} pose an attractive alternative to discretizing the domain using traditional body-fitted meshes. 

Solution methods for the Poisson equation on irregular domains have a long history in general numerical analysis, including the use of immersed methods. Most classical techniques are originally proposed based on second order discretizations, including the Shortley-Weller finite difference method \cite{Shortley1938, Jomaa2005}, Mayo's approach based on integral equations \cite{Mayo1984}, the original immersed interface method \cite{Leveque1994}, the cut-cell finite volume method \cite{johansen1998cartesian, Crockett2011}, and the ghost fluid method \cite{Fedkiw1999, Gibou2002}. Since these original seminal works, the respective methods have undergone substantial evolution and analysis, including their extensions to discretizations orders higher than second. Below we will briefly discuss advances in high order finite difference and finite volume methods, but we note briefly that immersed finite elements \cite{dePrenter2023}, discontinuous galerkin \cite{Saye2017}, and methods incorporating global extensions \cite{Stein2016} and integral equations \cite{Marques2017} are also active areas of research. 

In the context of finite difference methods, \cite{Ito2005} developed second and fourth order immersed interface stencils for handling Dirichlet and Neumann boundary conditions, and validation is shown in 2D. A 2D fourth order Poisson solver for Dirichlet boundary conditions was also developed in \cite{Linnick2005}, based on compact differences and an immersed interface method. In \cite{Gibou2005}, the authors developed a ghost cell method to achieve fourth order discretizations of the Laplace operator with Dirichlet boundary conditions in 2D domains. The Matched Interface and Boundary Method (MIBM) presented in  \cite{Zhou2006, Zhou2006a} was demonstrated up to sixth order accuracy on irregular 2D domains with immersed domain boundaries and interfaces. In \cite{Marques2011} the authors discretize and solve constant-coefficient Poisson interface problems to fourth order accuracy in 2D domains using compact difference schemes and a correct function method. Further, \cite{Zhu2016} considers a third order 2D Navier-Stokes solver based on compact differences with Dirichlet boundary conditions. Lastly, \cite{Hosseinverdi2018} presents a compact difference method for the Poisson equation with Dirichlet boundary conditions, achieving fourth order accuracy on 2D structured, graded grids. An extension to 3D by the same group was presented in \cite{Hosseinverdi2020}.

For cut-cell finite volume methods, fourth order elliptic discretizations have been developed for the Poisson equation with Dirichlet and Neumann conditions \cite{devendran2017fourth} and interface jump conditions \cite{Thacher2023AHigh}. These discretizations have been extended to more complex elliptic PDEs, including the unsteady Stokes equations on 2D or 3D domains \cite{OvertonKatz2022a} and to the 2D nonlinear shallow-shelf equations used to model ice sheet dynamics \cite{thacher2024high}.

In all these methods, the discretized problem leads to a sparse linear equation whose condition number increases at least quadratically with the inverse of the grid spacing. Solving a large-scale 3D system on fine meshes thus requires iterative solvers with appropriate conditioning strategies. Across the various reported 3D solvers in literature, different preconditioning and solution methods are used. In \cite{Linnick2005} the linear system was solved with a  multigrid method using an incomplete line LU (ILLU) decomposition smoother and an ad-hoc procedure for incorporating jumps. \cite{Zhou2006} uses the diagonal part of the discrete matrix as preconditioner for a BiCGStab method. In \cite{Hosseinverdi2020}, a low order discretization is solved using the modified strongly implicit method, and used to precondition a high order BiCGStab solver. In \cite{devendran2017fourth, OvertonKatz2022a} algebraic multigrid (AMG) or block-Jacobi preconditioners are used, and \cite{Devendran2014} presents a geometric multigrid (GMG) solver for a fourth order cut-cell discretization on adaptive grids in two dimensions. 

We further mention a number of relevant fast iterative solvers in the realm of second order immersed interface methods. In \cite{Adams:2002,Adams:2004}, the multigrid method is adapted to a classical second order immersed interface discretizations in 2D domains, by providing custom interpolation and restriction operators. Further, \cite{Chen:2008} provides a GMRES accelerated immersed multigrid solver for elliptic interface problems with second order discretizations in 2D domains, where the multigrid solver is made robust to complex geometries by adding multiple intersection points and custom interpolation operators. A uniform-grid alternative to the multigrid method is presented in \cite{Gillis2018}, using a low-rank decomposition of a second order immersed interface discretization based on the Sherman-Morrison-Woodbury (SMW) identity. The authors solve the global system using a Fast Fourier Transform with a Lattice Green's Functions kernel, while applying a GMRES-based solver for the low-rank problem. 

Based on the above literature, we note that there is a lack of methodologies providing both a high order immersed discretization of the Poisson equation, and an efficient solution technique that allows to solve large scale 3D problems. A possible reason is that the complexity of the discretization stencils in high order, 3D implementations poses challenges to directly assembling the matrix and applying efficient preconditioning strategies such as the algebraic multigrid method. Moreover, the stencil size of high order immersed methods typically involves some geometric constraints on the immersed geometry that can be resolved at a given resolution. These constraints may easily be violated on coarsened grids, further complicating the direct use of multigrid methods.  All of this is exacerbated when considering multiresolution grids, which are often required for efficient solution of practical 3D problems. On multiresolution grids, ghosts points or cells for differential operators are typically created on-the-fly using multidimensional interpolation stencils, which challenges methods that rely on explicitly assembled matrices. 

In this work we present a high order immersed interface method for the constant-coefficient Poisson equation with boundary and interface condition. The linear system arising from the high order immersed discretization is solved with a matrix-free GMRES approach, preconditioned using a novel geometric multigrid algorithm. The multigrid algorithm relies on a low order Shortley-Weller discretization adapted to Dirichlet, Neumann, and material interface conditions on the immersed geometry. With this combination our approach can handle high order immersed discretizations effectively preconditioned using a robust multigrid method that poses no additional constraints on the geometry. We demonstrate that our approach provides solution accuracies in complex domains up to sixth order for Dirichlet, Neumann, and interface conditions; compatibility with high order adaptive grid techniques; and parallel scalability for 3D problems up to \num{32768} compute cores. The resulting solutions presented here comprise, to the best of our knowledge, the largest high order solutions to elliptic problems solved using an immersed method to date.

The work is structured as follows. In section~\ref{sec:method} we first describe the immersed interface discretization based on earlier work on parabolic problems \cite{Gabbard2023, Gabbard2024}. Starting in section~\ref{subsec:linearsystem} we focus specifically on the linear system and analyze the discretization properties of our stencils for elliptic PDEs with various immersed boundary and interface conditions. In section~\ref{sec:solver}, we describe and analyze the fast iterative solver, with particular emphasis on the multigrid preconditioner. Lastly, in section~\ref{sec:results} we demonstrate the ability of our approach to generate high order and large scale 3D results with immersed domain boundary conditions, on both uniform and multiresolution grids.

\section{Discretization and numerical analysis}\label{sec:method}
In this section we will describe the proposed immersed interface discretization of the Poisson equation and perform an analysis of the truncation and discretization errors. For simplicity, all results in this section are presented for 2D domains and uniform grids, so that linear solutions can be found simply by forming a matrix and inverting it directly. The iterative solver is presented further below in section~\ref{sec:solver}.

\subsection{Discretization of immersed boundaries}\label{subsec:poisson_iim}

For immersed \textit{boundary} problems, we consider the Poisson equation posed in domain $\Omega^+$ with Dirichlet or Neumann conditions prescribed on the boundary, so that the PDE and boundary conditions are 
\begin{equation}\label{eq:poisson}
\begin{aligned}
    \nabla \cdot (\beta \nabla u) &= f \dbstext{in} \Omega^+, \\
    u &= \bar{u} \dbstext{on} \Gamma_D, \\
    \beta \partial_n u &= \bar{q} \dbstext{on} \Gamma_N.
\end{aligned}
\end{equation}
Here the prescribed boundary conditions and the source term $f(\xv)$ are assumed to be smooth functions, and the coefficient $\beta$ is assumed constant. The discretization of the Laplace operator follows the methodology outlined in \cite{Gabbard2023, Gabbard2024}, which is briefly summarized here. 

For interior points, the Laplacian operator $\nabla^2 = \sum_{i=1}^d \partial_i^2$ is discretized at each grid point using dimension splitting, with the standard central centered finite difference stencil used to discretize the second derivative $\partial_i^2$ along each coordinate axis. Near immersed boundaries an interpolation procedure incorporating boundary information is used to extend the solution, providing ghost values for the finite difference scheme. 

\begin{figure}[tb!]
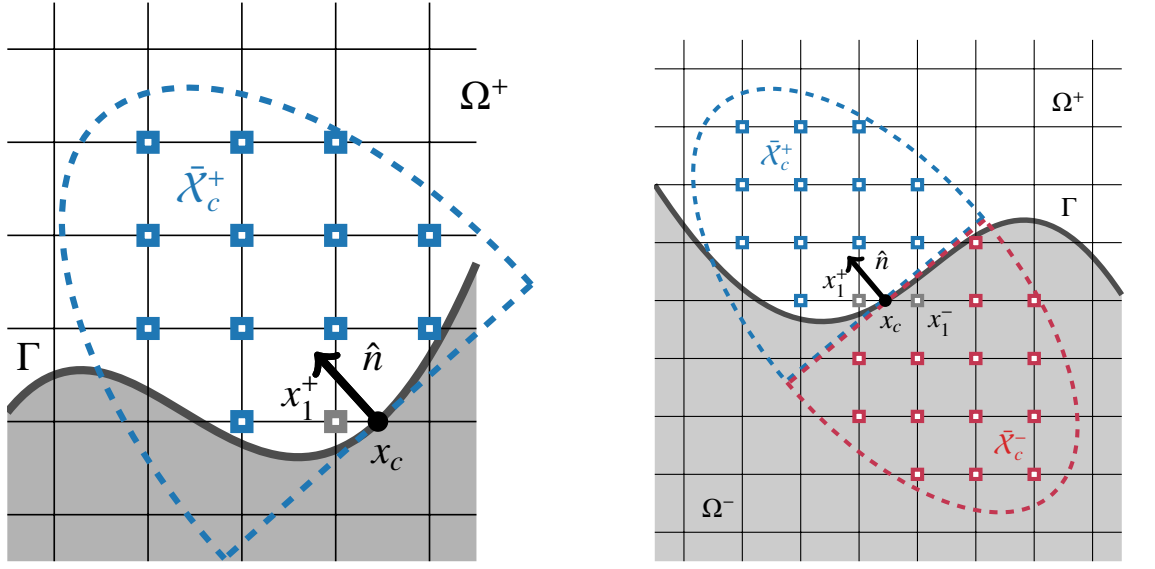

    \centering
    \begin{subfigure}{0.46\textwidth}
        \centering
        \resizebox{\textwidth}{!}{
            \input{tikzfigures/elliptical_diagram_onesided}
        }
        \caption{Multidimensional interpolant for immersed boundaries}
        \label{fig:iim_stencil_boundary}
    \end{subfigure}
    \hspace{0.05\textwidth}
    \begin{subfigure}{0.46\textwidth}
        \centering
        \resizebox{\textwidth}{!}{
            \input{tikzfigures/elliptical_diagram_twosided}
        }
        \caption{Multidimensional interpolants for immersed interfaces.}
    \label{fig:iim_stencil_interface}
    \end{subfigure}
    \caption{Each crossing between a grid line and the boundary ($x_c$)  is used to construct ghosts points for the affected grid points (light grey) using a multidimensional interpolant constructed from a half-elliptical region of grid points (half-ellipsoidal in 3D). For immersed boundaries (a), the interpolant is constructed using imposed boundary conditions (Dirichlet, Neumann). For immersed interfaces (b), polynomials from both sides are constructed using imposed jump conditions.}
    \label{fig:iim_stencil}
\end{figure}

The boundary treatments used in this work are based on the immersed interface method \cite{Leveque1994, Li2006, Wiegmann2000}, which maintains accuracy near immersed interfaces by accounting for known jumps in the solution and its derivatives. While this method was originally formulated using jump conditions derived from the continuous PDE, later reformulations \cite{Linnick2005, Gillis2018, Hosseinverdi2017} reduce the complexity of the method by using interpolation to estimate jumps in the solution. In many cases these reformulated methods are equivalent to an even simpler approach based on polynomial extrapolation \cite{Gabbard2022}, which we outline below. Following a convention from the immersed interface literature, we refer to the intersection between a 1D finite difference stencil and the surface $\Gamma$ as a \textit{control point}, denoted $\xv_c$. The evaluation point for this stencil is referred to as an \textit{affected point}, since the discretization is affected by the presence of the surface. 

In our approach, each control point $\xv_c$ on the immersed domain boundary is associated with a set of interpolation points 
$\mathcal{X}_c^+ \subset \Omega^+$, 
and with a polynomial $p_c(\xv)$ that approximately interpolates the domain values $\qty{u(\xv_i) \mid \xv_i \in \mathcal{X}_c^+}$ in a least squares sense. Each 1D finite difference stencil that intersects the boundary at $\xv_c$ is applied not directly to the function $u(\xv)$, but to the extended function
\begin{equation}\label{ch2:eq:extended_function}
    u_c(\xv) = \begin{cases} u(\xv), & x \in \Omega \\ p_c(\xv), & x \notin \Omega \end{cases}.
\end{equation}

As defined in \cite{Gabbard2024}, the set of interpolation points in the least squares domain  $\mathcal{X}_c^+$ includes the control point, excludes the closest grid point, and includes all other grid points that fall within a half-elliptical region centered on the boundary (half-ellipsoidal in 3D) (Fig.~\ref{fig:iim_stencil_boundary}). These interpolants are suitable for all boundary conditions and any smooth geometry satisfying a well-defined curvature constraint \cite{Gabbard2024}. In this work, for any given polynomial order we choose the major and minor axes of the region so that we can guarantee the existence of $p_c(\xv)$ on a grid with spacing $\Delta x$ as long as the immersed surface satisfies
\begin{equation}
\label{eq:curvature_condition}
|\kappa \Delta x| < 1/4,
\end{equation}
where $\kappa$ is the maximum scalar curvature in 2D, or the maximum of the two principal curvatures in 3D. 

Each of the 1D finite difference stencils which intersect the boundary at $\xv_c$ can require up to $w$ ghost values at grid points $\{\xv_g\}$ that fall outside of the domain. These ghost values are obtained by evaluating $p_c(\xv)$ through $w$ separate stencil operations with coefficients $\{s_c^g, s^g_\alpha\}_{\alpha = 1}^n$ \cite{Gabbard2024}, so that 
\begin{equation}
    p_c(\xv_g) = s^g_c u(\xv_c) + \sum_{\alpha = 1}^n s^g_\alpha u(\xv_\alpha)
\end{equation}
at each point $\xv_g$ requiring a ghost value. For points with Neumann boundary conditions $u(\xv_c)$ is not directly available, but it can be approximated based on the boundary condition $\bar{q}(\xv_c)$ and nearby solution values. To clarify,  let $\{s_c, s_i\}_{i = 1}^n$ be a set of stencil coefficients that approximate the normal derivative of $p_c$ at $\xv_c$, so that
\begin{equation}\label{ch2:eq:normal_grad}
    \beta \partial_n u(\xv_c) = s_c u(\xv_c) + \sum_{i = 1}^{n} s_i u(\xv_i) + \order{\Delta x^{k - 1}}.
\end{equation}
When a Neumann condition is prescribed, Eq.~\ref{ch2:eq:normal_grad} can be inverted to give
\begin{equation}\label{ch2:eq:neumann}
    u(\xv_c) = \frac{1}{s_c} \qty(q(\xv_c) - \sum_{i = 1}^{N-1} s_i u(\xv_i)) + \order{\Delta x^{k}}.
\end{equation}
This requires one additional set of stencil coefficients which evaluate the normal derivative $\partial_n p(\xv_c)$ on the boundary. 

Computing the full set of stencil coefficients for the multivariant interpolant requires a single QR factorization at each control point, followed by one triangular solve and one matrix-vector product per stencil (a total of $w$ for Dirichlet boundary conditions and $w + 1$ for Neumann boundary conditions).

\subsection{Discretization of immersed interfaces}
For problems with immersed interfaces, we consider the Poisson equation defined over two subdomains $\Omega^+$ and $\Omega^-$ with piecewise constant coefficients $\beta^+$ and $\beta^-$:
\begin{equation}\label{eq:piecewise_beta}
    \beta(\xv) = \begin{cases} \beta^+, & \xv \in \Omega^+ \\ \beta^-, & \xv \in \Omega^- \end{cases}.
\end{equation}
Jump conditions on the solution and its flux are prescribed on the interface between the subdomains, $\Gamma_M$, so that the full problem becomes
\begin{equation}\label{eq:poisson_jumps}
\begin{aligned}
    \nabla \cdot (\beta\nabla u) &= f \dbstext{in} \Omega, \\
    [u] &= j_0(\sv) \dbstext{on} \Gamma_M, \\
    [\beta \partial_n u] &= j_1(\sv) \dbstext{on} \Gamma_M.
\end{aligned}
\end{equation}
The source term $f(\xv)$ is also allowed to be discontinuous across the material interface, though this has minimal effect on the solution procedure. Equations~\ref{eq:piecewise_beta} and \ref{eq:poisson_jumps} often occur in simulations containing two materials with distinct properties, such as steady-state heat conduction in materials that have different thermal conductivities or electrostatics with materials that have differing permittivities. Though the classical IIM was developed exactly for this setting \cite{Leveque1994, Li1997}, our approach builds upon our immersed boundary treatment described above to enable (1) robust extension to concave geometries in 2D and 3D, (2) straightforward extension to high order, and (3) consistent implementation across all boundary/interface conditions. Though not relevant here, our approach has the additional flexibility of allowing \textit{weighted} least squares polynomial fits, which significantly improves the stability of high order discretizations of hyperbolic PDEs in 2D and 3D \cite{devendran2017fourth, Gabbard2023, Gabbard2024}. 

For piecewise constant $\beta(\xv)$, the operator $\nabla \cdot (\beta \nabla u)$ reduces to $\beta \nabla^2 u$ away from $\Gamma_M$, so that the standard dimension-split discretization of interior points remains valid. To discretize the jump boundary conditions, the boundary values $u^-(\xv_c)$ and $u^+(\xv_c)$ from either side of the interface are computed with the aid of two sets of stencils coefficients \cite{Gabbard2023,Gabbard2024}, associated with half-elliptical regions as shown in Fig~\ref{fig:iim_stencil_interface}. The first set $\{s^+_c, s^+_i\}$ maps the boundary value $u^+(\xv_c)$ and solution values from $\Omega^+$ to the normal derivative $\partial_n u^+(\xv_c)$, while the second set $\{s_c^-, \, s_i^-\}$ is designed analogously to map solution values from $\Omega^-$ to the normal derivative $\partial_n u^-(\xv_c)$. The boundary values $u^\pm(\xv_c)$ can then be determined from the discretized jump conditions $u^+(\xv_c) - u^-(\xv_c) = j_0(\xv_c)$ and
\begin{equation}
\label{ch2:eq:boundary_system_scalar}
        \beta^+\qty(s^+_c u^+(\xv_c) + \sum_{i=1}^{n^+} s^+_i u(\xv_i^+)) - \beta^-\qty(s^-_c u^-(\xv_c) + \sum_{i=1}^{n^-} s^-_i u(\xv_i^-)) = j_1(\xv_c).
\end{equation}
For $j_0(\xv_c) = j_1(\xv_c) = 0$, the closed-form solution to Eq.~\ref{ch2:eq:boundary_system_scalar} is 
\begin{equation}
   u^-(\xv_c) = u^+(\xv_c) = \bar{u} \dbstext{with} \bar{u} = -\frac{
        \beta^+ \sum_{i=1}^{n^+} s^+_i u(\xv_i^+) - 
        \beta^- \sum_{i=1}^{n-} s^-_i u(\xv_i^-)
    }{\beta^+ s_c^+ - \beta^- s_c^-}.
\end{equation}
When jumps are present, the boundary values instead are given separately by
\begin{equation}\label{ch2:eq:explicit_material_interface}
    u^+(\xv_c) = \bar{u} + \frac{j_1 - \beta^- s_c^- j_0(\xv)}{\beta^+ s_c^+ - \beta^- s_c^-}, \quad u^-(\xv_c) = \bar{u} + \frac{j_1 - \beta^+ s_c^+ j_0(\xv)}{\beta^+ s_c^+ - \beta^- s_c^-}.
\end{equation}
The computation of boundary stencils for each point on a material interface requires one QR factorization, $w + 1$ matrix-vector products, and $w + 1$ triangular solves for each side of the interface, exactly twice the number of operations needed per point compared with Neumann domain boundaries. Once determined, the boundary values $u^\pm(\xv_c)$ can be used in stencil operations on either side of the interface. We note that Equation~\ref{ch2:eq:boundary_system_scalar} treats both sides of the interface on an even footing. Further, as the ratio $\beta^- / \beta^+$ tends to zero, Eq.~\ref{ch2:eq:explicit_material_interface} approaches the well-behaved Neumann boundary treatment given the previous section. This indicates that the boundary treatment is robust to large jumps in coefficients and does not exhibit singular behavior when $\beta^-/\beta^+$ tends to zero or infinity.

\subsection{Linear system}
\label{subsec:linearsystem}
The interior discretization, treatment of domain boundaries, and treatment of material interfaces completes the full discretization of the Poisson equation. The methodology outlined above allows for a range of high order Poisson discretizations, each defined by the choice of interior stencil and the order of accuracy of the boundary interpolants $p_c(\xv)$. For brevity, we refer to a discretization with an $n$-th order interior stencil and $k$-th order boundary interpolants as an $(n, k)$  discretization, and focus on high order interior discretizations with $n = 4$ or $n = 6$.

After discretization, the Poisson equation can be expressed in the form of a system of linear equations relating the solution values at all grid points in $\Omega$ to the known boundary data and source data. After all stencil coefficients have been computed, both the interior and boundary discretizations are assembled into a large sparse linear system of the form
\begin{equation}\label{eq:discrete_poisson}
    \vb{L}_{\Omega} \vb{u}_{\Omega} + \vb{L}_{\Gamma} \vb{u}_{\Gamma} = \vb{f}_{\Omega}
\end{equation}
Here $\vb{u}_{\Omega}$ contains the unknown solution at the grid points, $\vb{u}_{\Gamma}$ contains the known boundary or interface data ($\bar{u}$, $\bar{q}$, or $j_0$ and $j_1$) at the control points, $\vb{f}_{\Omega}$ contains the known source term $f(\xv)$ at the grid points, and the entries of the matrices $(\vb{L}_{\Omega},\,\vb{L}_{\Gamma})$ are derived from the interior and boundary/interface stencil coefficients. Generally, $\vb{L}_{\Omega}$ is not symmetric. 

For singular problems the solution to the continuous PDE is defined only up to an additive constant, and the source term $f(\xv)$ must obey an integrability condition. For purely Neumann boundary conditions, for instance, the integrability condition is posed as $\int_{\Omega} f(\xv) \dd{\vb{x}} = - \int_{\Gamma_N} \bar{q}(\sv) \dd{\vb{s}}$. When applied to these problems, the system matrix $\vb{L}_{\Omega}$ is singular and satisfies $\vb{L}_{\Omega} \mathds{1} = 0$, where $\mathds{1}$ is a vector of ones. As a result, the discrete solution is also defined only up to an additive constant. Following \cite[Chapter~5.6.4]{Trottenberg2000} the linear system can be augmented with an additional unknown to explicitly account for this constant, leading to the augmented system
\begin{equation}
    \mqty[\vb{L}_{\Omega} & \mathds{1} \\ \mathds{1}^T & 0] \mqty[\vb{u}_{\Omega} \\ \alpha] = \mqty[\vb{f}_{\Omega} - \vb{L}_{\Gamma} \vb{u}_{\Gamma} \\ \gamma].
\end{equation}
The scalar $\gamma$ controls the additive constant directly via $\mathds{1}^T \vb{u}_{\Omega} = \gamma$, while $\alpha$ is an unknown shift determined so that $\vb{f}_{\Omega} - \vb{L}_{\Gamma} \vb{u}_{\Gamma} + \alpha \mathds{1}$ is in the column space of $\vb{L}_{\Omega}$. This latter condition is the discrete analog of the continuous integrability condition. 

Below, for convergence tests of singular problems, the exact and numerical solution vectors are adjusted to have zero average before they are compared, so that reported errors in the solution do not reflect the additive constant.

\subsection{Truncation error and order of accuracy}\label{subsec:poisson_error}
\begin{figure}[tb!]
    \centering
    \begin{subfigure}{0.49\textwidth}
        \centering
        \resizebox{\textwidth}{!}{\input{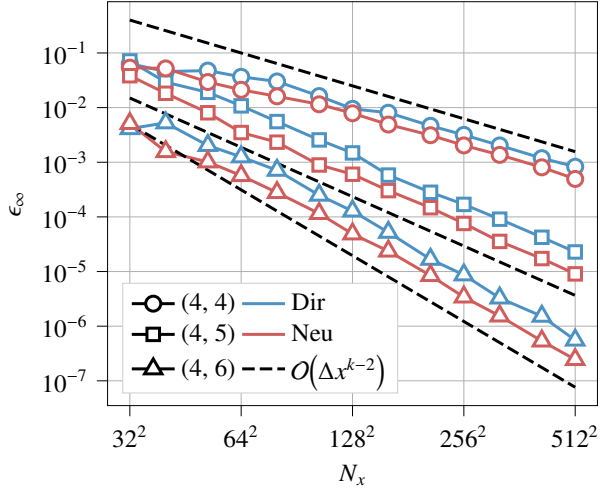}}
        \caption{Truncation error of (4, k) schemes}
        \label{subfig:poisson_truncation_four}
    \end{subfigure}
    \begin{subfigure}{0.49\textwidth}
        \centering
        \resizebox{\textwidth}{!}{\input{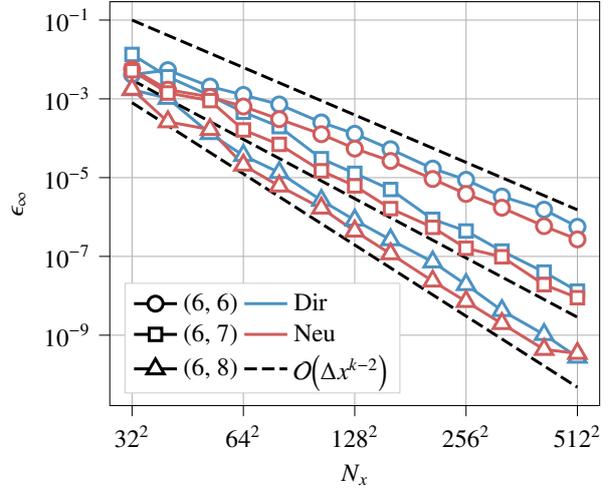}}
        \caption{Truncation error of (6, k) schemes}
        \label{subfig:poisson_truncation_six}
    \end{subfigure}
    \caption{Truncation error of Poisson discretizations with (a) fourth order and  (b) sixth order interior schemes as a function of the grid spacing $\Delta x$ for both Dirichlet and Neumann boundary conditions. In all cases the truncation error of each $(n, k)$  discretization is of order $\Delta x^{k - 2}$, indicated by dashed lines, with the Neumann cases showing a slightly lower error magnitude than their Dirichlet counterparts.}
    \label{fig:poisson_trunction}
\end{figure}

For an $(n, k)$  discretization of the Poisson equation, the truncation error at interior points is that of the interior difference stencil, which is of order $\Delta x^n$. The truncation error at the affected points is of order of order $\Delta x^{\mathrm{min}(n, k - 2)}$, which depends on both interior stencil and the accuracy of polynomial interpolants $p_c(\xv)$. Thus to maintain a truncation error of order $n$ throughout the domain, it is necessary to use interpolants of order $k = n + 2$. To illustrate this, consider a test case with the manufactured solution
\begin{equation}\label{eq:poisson_manufactured}
    u(\xv) = \sin(4\pi x_1) \sin(2\pi x_2).
\end{equation}
Let the domain $\Omega$ be the exterior of the star-shaped test geometry defined in \cite{Gabbard2024}, with center $\xv_0 = [0.501, 0.502]$, average radius $r_0 = 0.28$, and perturbation radius $\tilde{r} = 0.025$. Figure~\ref{fig:poisson_trunction} plots the $L_\infty$ norm of the truncation error that results from applying ($n$, $n$), ($n$, $n + 1$), and ($n$, $n + 2$)  Poisson discretizations to this test case with $n = 4, 6$. The test is repeated over a range of spatial resolutions, characterized by the number of grid points $N_x$ along each axis, and with both Dirichlet and Neumann boundary conditions. In all cases each $(n, k)$ discretization has a truncation error of order $k - 2$.

\begin{figure}[tb!]
    \centering
    \begin{subfigure}{0.49\textwidth}
        \centering
        \resizebox{\textwidth}{!}{\input{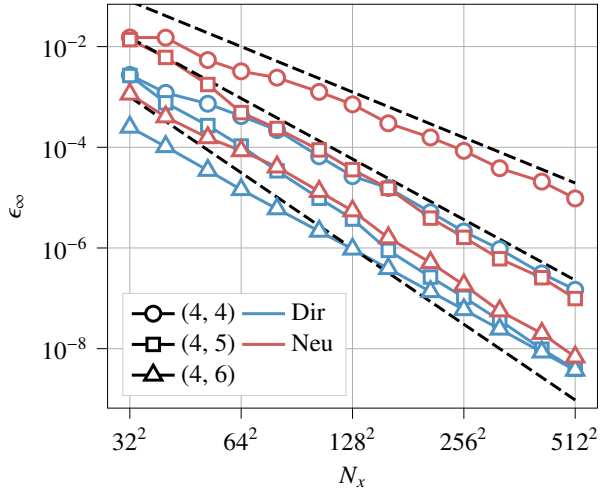}}
        \caption{Solution error of (4, k) schemes}
        \label{subfig:poisson_solution_four}
    \end{subfigure}
    \begin{subfigure}{0.49\textwidth}
        \centering
        \resizebox{\textwidth}{!}{\input{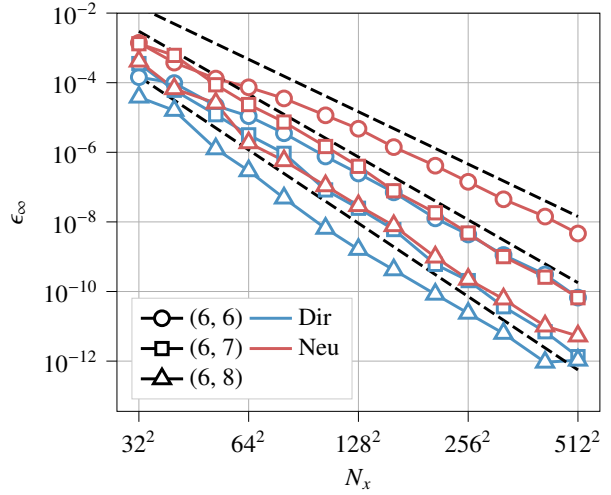}}
        \caption{Solution error for (6, k) schemes}
        \label{subfig:poisson_solution_six}
    \end{subfigure}
    \caption{Solution error (4, k) and (6, k) Poisson discretizations with both Dirichlet and Neumann boundary conditions. The black dashed lines in each plot indicate slopes of $\order{\Delta x^{n - 1}}$, $\order{\Delta x^{n}}$, and $\order{\Delta x^{n+1}}$, from top to bottom.}
    \label{fig:poisson_solution}
\end{figure}

The convergence of the discrete solution $u(\xv_i)$ depends on both the truncation error and on the nature of the prescribed boundary conditions. Generally an $(n, k)$ discretization cannot converge faster than $\Delta x^n$, but may still converge faster than the boundary truncation error of order $\Delta x^{\mathrm{min}(n, k - 2)}$. This phenomenon was previously observed and explained for cut-cell finite volume methods in \cite{johansen1998cartesian}, and the results observed here for finite difference methods follow a similar pattern. Figures~\ref{subfig:poisson_solution_four} and \ref{subfig:poisson_solution_six} plot the $L_\infty$ error in the discrete solution obtained using the same manufactured solution as above with either a fourth or sixth order interior discretization. For interior order $n$ and Dirichlet boundary conditions, the $(n, n)$, $(n, n+1)$, and $(n ,n+2)$ discretizations all display convergence of order $n$, with a prefactor that decreases for higher order boundary interpolants. In regimes where the error associated with the boundaries dominates, the convergence of the $(n, n+1)$, and $(n, n+2)$ may appear to be even faster; however, when carried to a fine enough resolution the interior error of order $n$ ultimately dominates. For Neumann boundary conditions, the $(n, n)$ discretization achieves order $n - 1$ convergence, while both the $(n, n + 1)$, and $(n, n + 2)$ discretizations converge slightly faster than order $n$ (though clearly slower than order $n + 1$). The slight superconvergence can again be attributed to boundary error which converges at a faster rate than the interior error.

Altogether, these results indicate that an $(n, n)$ discretizations will only achieve $n$-th order accuracy with Dirichlet boundary conditions. An $(n, n + 1)$ discretization achieves $n$-th order accuracy with both Dirichlet or Neumann boundary conditions, and leads to lower overall error magnitudes even with Dirichlet boundary conditions. While an $(n, n + 2)$ discretization can reduce the magnitude of the error even further, it does not bring any further increase in convergence order.

In addition to the convergence of the discrete solution, the convergence of boundary quantities as well as the convergence of the gradient of the solution are important in many applications. Figure~\ref{fig:poisson_boundary} plots the $L_\infty$ error in the non-prescribed boundary quantity for simulations with (4, 5) and (6, 7)  discretizations. For Dirichlet discretizations this is the error in the boundary normal gradient, calculated by applying a boundary stencil operation at each control point. For Neumann discretizations this the error in the solution on the boundary, calculated as in Eq.~\ref{ch2:eq:neumann}. In all cases the convergence rate is of order $n$ for an ($n$, $n + 1$) discretization, demonstrating the ability of high order discretizations to accurately predict quantities defined on an immersed surface. Moreover, Figure~\ref{fig:poisson_gradient} plots the $L_\infty$ error of the gradient of the solution for simulations with (4, 5) and (6, 7)  discretizations. Here each component of the gradient is computed using a centered finite-difference stencil for the first derivative. Near the boundary, the same immersed interface treatment as described for the Laplacian in section~\ref{subsec:poisson_iim} is used. The prescribed boundary conditions (Dirichlet or Neumann) are incorporated in the immersed interface stencil. For each case $(n,k)$ in Figure~\ref{fig:poisson_gradient}, the gradient is computed using the same $(n,k)$ discretization orders as the solution itself. The results show that the solution gradient converges at the same order as the solution itself, exhibiting superconvergence as discussed in, amongst others, \cite{Gibou2002,Gallinato2015, Yoon2015}.

\begin{figure}[tb!]
\begin{minipage}[t]{0.48\textwidth}
        \centering
        \resizebox{\textwidth}{!}{\input{results/stencil/poisson_boundary_2D}}
        \caption{$L_\infty$ norm convergences of the unknown boundary quantity for (4, 5) and (6, 7)  discretizations with Dirichlet or Neumann boundary conditions. For Dirichlet discretizations the unknown is the boundary normal gradient, while for Neumann discretizations it is the boundary solution values.}
        \label{fig:poisson_boundary}
\end{minipage}
\hfill
\begin{minipage}[t]{0.48\textwidth}
\centering
\resizebox{\textwidth}{!}{\input{results/stencil/poisson_gradient_2D}}
        \caption{$L_\infty$ norm convergences of the solution gradient for (4, 5) and (6, 7)  discretizations with Dirichlet or Neumann boundary conditions. The solution gradient is computed using the same immersed interface discretization as used for the Laplacian, using a centered finite difference stencil.}
        \label{fig:poisson_gradient}
\end{minipage}        
\end{figure}

\begin{figure}[tb!]
    \centering
    \begin{subfigure}{0.45\textwidth}
        \centering
        \includegraphics[width=\textwidth]{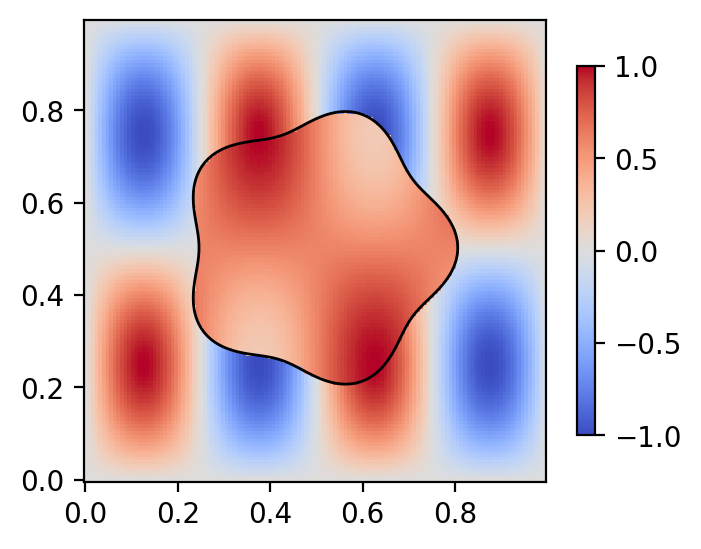}
        \caption{Manufactured solution}
        \label{subfig:material_setup}
    \end{subfigure}
    \begin{subfigure}{0.45\textwidth}
        \centering
        \resizebox{\textwidth}{!}{\input{results/stencil/material_combined_ratio_2}}
        \caption{Solution and boundary error}
        \label{subfig:material_error}
    \end{subfigure}
    \caption{Results for the Poisson equation with material interfaces. (a) the geometry and manufactured solution for this test case. (b) $L_\infty$ norm convergences of solution and the boundary normal gradient $\partial_n u^{+}(\sv)$ with $\beta^- / \beta^+ = 2$. The (4, 5) and (6, 7)  discretizations achieve fourth and sixth order convergence respectively for both quantities. }
    \label{fig:material_convergence}
\end{figure}

To test the material interface treatment, consider a test case with domains $\Omega^+$ and $\Omega^-$ occupying the exterior and interior of the same star-shaped body. Each region has coefficient $\beta^- = 1$ and $\beta^+ = 1/2$, with a manufactured solution
\begin{equation}\label{eq:manufactured_interface}
    u(\xv) = \begin{cases}
        0.6 + 0.4\sin(4\pi x_1) \sin(2\pi x_2), & \xv \in \Omega^+ \\
        \sin(4\pi x_1) \sin(2\pi x_2), & \xv \in \Omega^-
    \end{cases}.
\end{equation}
This setup is illustrated in Fig.~\ref{subfig:material_setup}. The $L_\infty$ error in the discrete solution is given in Fig.~\ref{subfig:material_error}, demonstrating fourth and sixth order spatial convergence for (4, 5) and (6, 7)  discretizations with material interfaces. Also shown is the $L_\infty$ error of the interface gradient $\partial_n u^+(\sv)$ for $\sv \in \Gamma^M$, which indicates that the solution and its gradient also converge at order $n$ on the material interface.

\subsection{Spectrum of the discretized Laplacian}\label{subsec:poisson_spectrum}
When the domain $\Omega$ is periodic and contains no immersed boundary, the spectrum of the spatial discretization can be completely characterized by a von Neumann analysis. The discrete operator consists only of a dimension-split finite difference scheme based on a standard symmetric second derivative finite difference stencil. In $d$ dimensions, this operator has eigenfunctions of the form $u(\xv) = \exp(\kv \cdot \xv)$ and corresponding eigenvalues $\lambda(\vb{k}) = -\frac{1}{\Delta x^2}\sum_{i = 1}^d \sigma(k_i \Delta x)$, where
\begin{equation}
    \sigma(k) = 2 \sum_{j = 1}^w a_j (1 - \cos(jk))
\end{equation}
and $\{a_j\}_{j = -w}^{w}$ are the coefficients of the one-dimensional finite difference stencil. Note that the function $\sigma(\cdot)$ is real, even, and $2\pi$-periodic, and for all stencils considered here $\sigma(k) > 0$ except at $k = 0$. As a result, the spectrum of the discrete operator consists of real eigenvalues satisfying
\begin{equation}
    \Delta x^2 \lambda(\kv) \in [-d \sigma_{\mathrm{max}}, 0], \dbstext{where} \sigma_{\mathrm{max}} = \max_{0 \le k < \pi} \sigma(k).
\end{equation}
The value of $\sigma_{\mathrm{max}}$ for each centered second derivative stencil is provided in Table~\ref{tab:sigma_max}.

\begin{figure}[tb!]
    \centering
    \begin{subfigure}{0.49\textwidth}
    \centering
        \includegraphics[width=\textwidth]{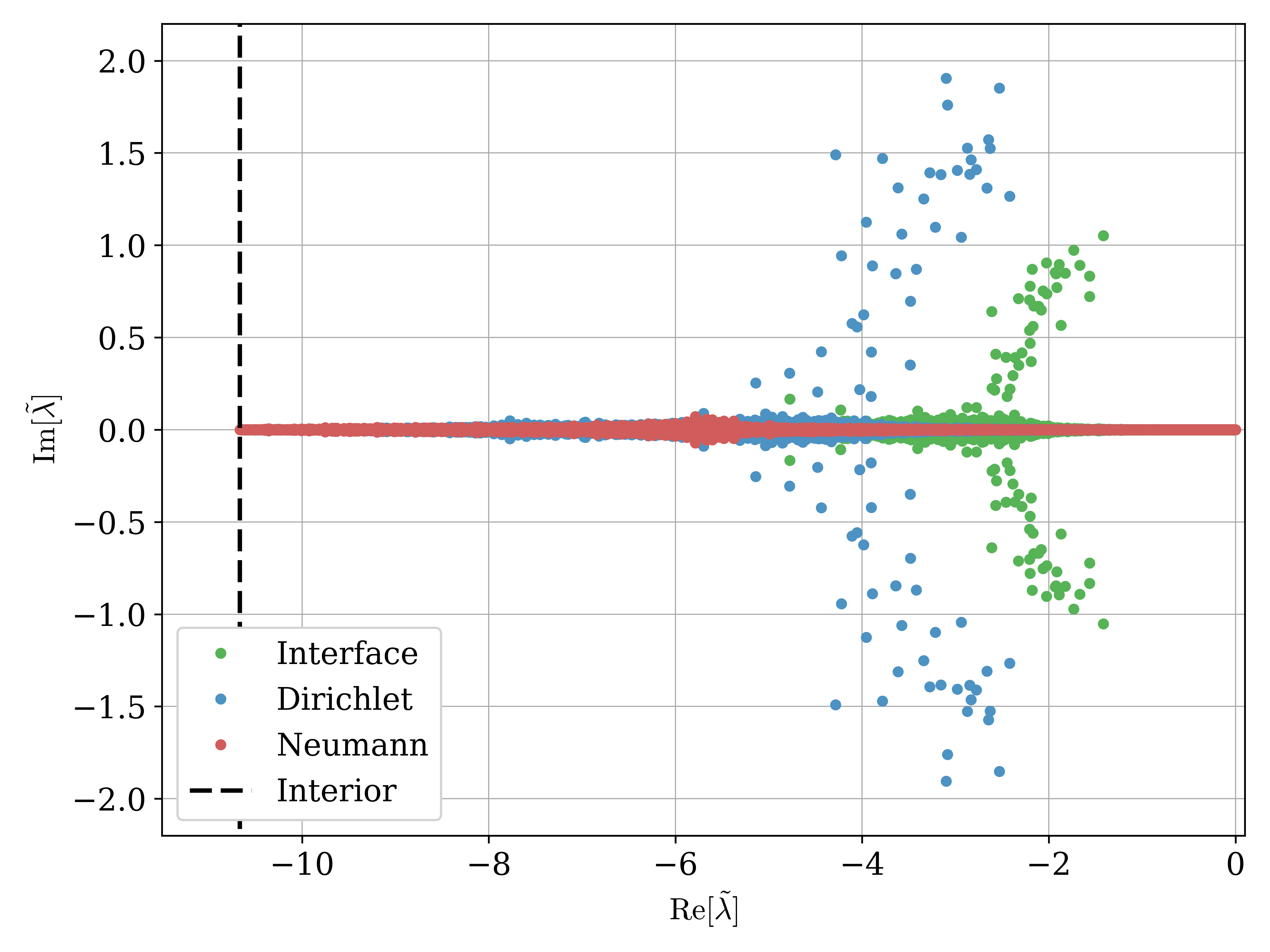}
        \caption{(Spectrum of the (4, 5) operator}
        \label{subfig:poisson_spectrum_four}
    \end{subfigure}
    \begin{subfigure}{0.49\textwidth}
        \centering
        \includegraphics[width=\textwidth]{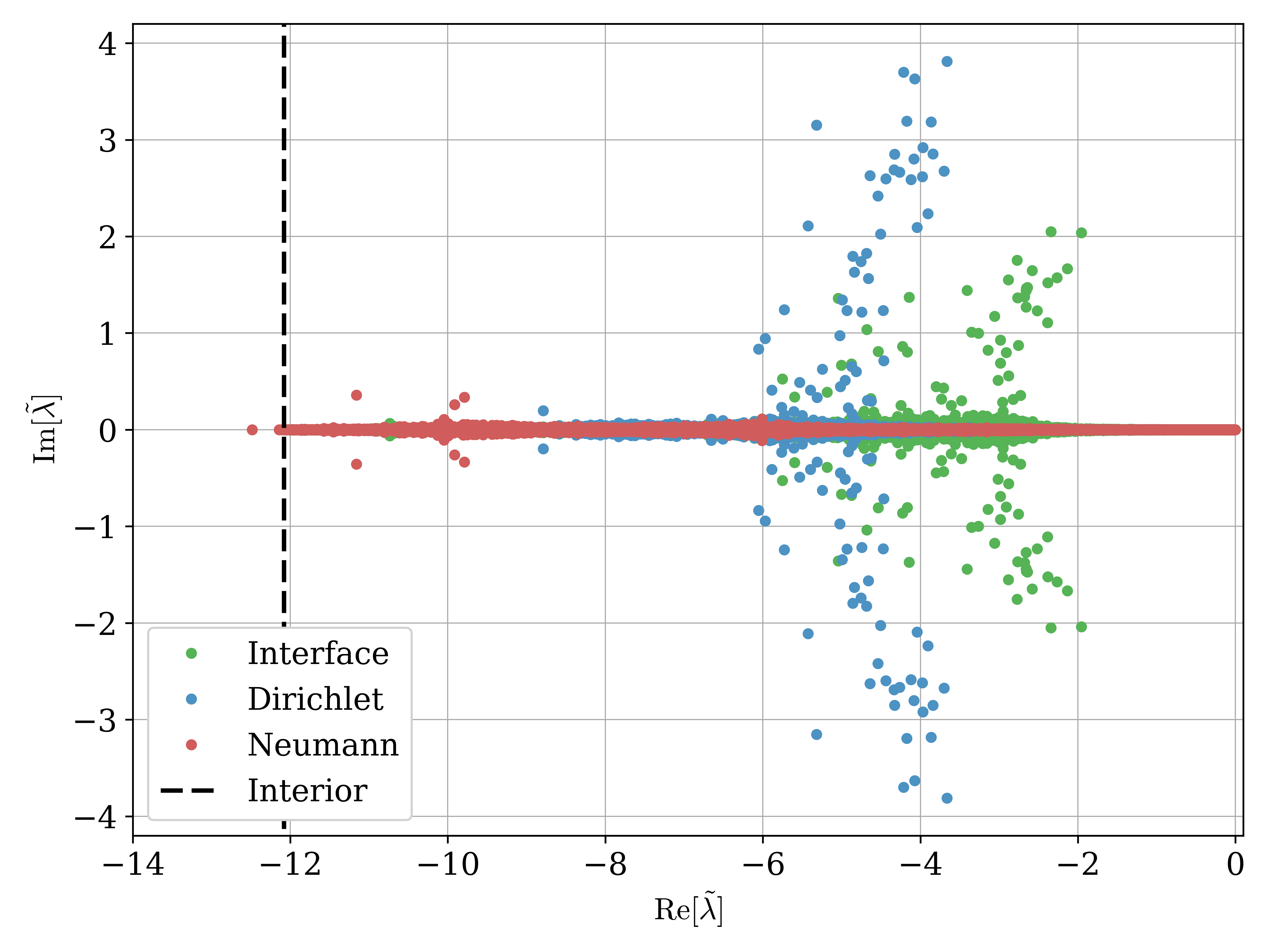}
        \caption{Spectrum of the (6, 7) operator}
        \label{subfig:poisson_spectrum_six}
    \end{subfigure}
    \caption{Spectrum of the (4, 5) and (6, 7) Poisson discretizations with Dirichlet (blue), Neumann (red), or interface jump (green) conditions on an immersed interface at $N_x = 64$. In all cases the eigenvalues satisfy $\Re \tilde{\lambda} \le 0$, and the eigenvalue with the most negative real part agrees closely with a von-Neumann analysis of the interior discretization. For both fourth and sixth order the Dirichlet case produces eigenvalues with larger imaginary components, while the corresponding Neumann case has eigenvalues that remain close to the negative real axis. The interface jump case spectra  are consistent with those of the Neumann discretization.}
    \label{fig:poisson_spectrum}
\end{figure}

When immersed boundaries are present, the system matrix $\vb{L}_{\Omega}$ is not symmetric, and as a result its spectrum is complex valued. Figure~\ref{fig:poisson_spectrum} plots the scaled eigenvalues $\tilde{\lambda} = \Delta x^2 \lambda$ for both (4, 5) and (6, 7) discretizations, using the 2D test geometry from the previous section and linear spatial resolution $N_x = 64$. For the interface jump conditions, the eigenvalues are further scaled by $\beta_{\max} = \max(\beta^-, \beta^+)$, so that $\tilde{\lambda} = \Delta x^2 \lambda/\beta_{\max}$. In each case the plot also includes a vertical line at $-2 \sigma_{\mathrm{max}}$ that indicates the most negative eigenvalue of the interior scheme. In all cases the eigenvalues satisfy $\Re \tilde{\lambda} \le 0$, and for all but the sixth order Neumann discretization the eigenvalues also satisfy $-d \sigma_{\mathrm{max}} \le \Re \tilde{\lambda}$. For the sixth order Neumann discretization the most negative eigenvalue exceeds that of the interior discretization, but only by 3.3\%. Taken together, these results indicate that the real part of the spectrum of $\vb{L}_{\Omega}$ is largely unaffected by the addition of immersed domain boundaries. Consequently, the discretization does not suffer from a ``small cell'' issue, which would result in extremely large negative eigenvalues. 

\begin{table}[htb!]
\caption{\label{tab:sigma_max} Maximum eigenvalue magnitude for centered second derivative stencils.}
\centering
\begin{tabular}{ c c c c } 
\toprule
Order & 2 & 4 & 6 \\
\midrule
$\sigma_\mathrm{max}$ & $4$ & $16/3 = 5.3\bar{3}$ & $272 / 45 = 6.0\bar{4}$\\
\bottomrule
\end{tabular}
\end{table}

The spectrum of $\vb{L}_{\Omega}$ is especially relevant when iterative linear solvers are applied to the discrete linear system, as discussed in section~\ref{sec:solver}. With immersed boundaries, the non-symmetric system matrix precludes the use of a preconditioned conjugate gradient method, which is often the standard for symmetric positive-definite Poisson discretizations produced with finite elements. Among the popular Krylov subspace methods, the lack of symmetry limits the available solvers to preconditioned GMRES or BiCGSTAB. With or without immersed boundaries, the eigenvalues of $\vb{L}_{\Omega}$ have magnitudes which range from $\inlineorder{1}$ to $\inlineorder{1/\Delta x^2}$, indicating a high condition number of $\kappa \sim 1/\Delta x^2$. Consequently, preconditioning is required to obtain a solution in a reasonable number of iterations. 

\section{Fast iterative solver}\label{sec:solver}
Here we will describe the components of the proposed fast iterative solver for immersed discretizations of the Poisson equation in two and three dimensions. Given the need for preconditioning, we start by explaining the multigrid preconditioners for domain boundary conditions (Dirichlet and Neumann) followed by the multigrid preconditioner for interface conditions (jump conditions). 

Before continuing, we note that in our 3D multiresolution implementation the structure and coefficients of the matrix $\vb{L}_{\Omega}$ depend on the particular treatment of resolution jumps in the domain. The multiresolution results in this work are based on \textsc{murphy} \cite{Gillis2022}, a software framework that relies on a wavelet-based multiresolution analysis, which makes the matrix $\vb{L}_{\Omega}$ difficult to assemble explicitly. Our approach is therefore guided by the need for a matrix-free implementation of the operators $\vb{L}_{\Omega}$ and $\vb{L}_{\Gamma}$. While this poses restrictions on the algorithms that can be used, it enables straightforward coupling with arbitrary multiresolution adaptation strategies and significantly simplifies the high order immersed interface treatment. 

\subsection{Shortley-Weller discretizations}
For our high order immersed interface stencil, which is defined under the assumption that %
geometric features are well resolved (Eq~\eqref{eq:curvature_condition}), successive coarsening will inevitably lead to ill-posed systems for formulating the multivariate polynomial interpolant. This motivates the use of a lower-order preconditioner that allows coarsening to arbitrary resolutions while still accounting for the boundary or interface conditions on the immersed surface.

An existing approach that offers this geometric flexibility is the Shortley-Weller (SW) discretization of the Laplacian \cite{Shortley1938}. The SW stencil coincides with a standard, 5 point (in 2D) or 7 point (in 3D) second-order finite difference stencil away from the boundary. Whenever the stencil crosses an interface, it is adapted to incorporate the boundary condition instead of the neighboring grid value. When applied to the Poisson equation with Dirichlet boundary conditions, this discretization yields solutions that are globally second order accurate in the $L_\infty$ norm and third order accurate at points near the boundary \cite{Matsunaga2000}. Because the stencil incorporates information from nearest neighbors only, it can be successfully applied to arbitrarily coarse geometries, and there are already successful multigrid methods \cite{Trottenberg2000} and preconditioned high order immersed methods \cite{Hosseinverdi2020} based on the SW discretization available in the literature. Here we build on these methods by incorporating low-order discretizations of Neumann and jump conditions, and by tuning the SW boundary discretization to better match the spectrum of our high order methods.

We introduce the Shortley-Weller discretization in 1D, and extend it to 2D and 3D domains in an entirely dimension split way. Consider a 1D uniform grid with grid points $\{x_i\}_{i \in \mathbb{Z}}$ and grid spacing $\Delta x$, and let $\{u_i\}$ be the values of a smooth function $u(x)$ at the grid points. A grid point $x_i$ is considered to be irregular if either $x_{i-1}$ or $x_{i+1}$ lies outside the problem domain, and regular otherwise. At regular points the Laplacian is approximated by the standard stencil $\partial_x^2 u \approx (u_{i-1} - 2u_i + u_{i + 1})/\Delta x^2$. For an irregular point $x_i$, we distinguish between three cases: (a) $x_{i-1}$ lies outside the domain, (b) $x_{i + 1}$ lies outside the domain, and $(c)$ both $x_{i - 1}$ and $x_{i+1}$ lie outside the domain. For cases (a) and (c) let $x_w^- \in [x_{i-1},\, x_i]$ denote the location of the domain boundary and let $\Delta x^- = x_i - x_w^-$ be the distance to the boundary point. Likewise for case (b) and (c) let $x_w^+ \in [x_i,\,x_{i+1}]$ be the boundary point and let $\Delta x^+ = x_w^+ - x_i$. In terms of the fractional grid spacing $\psi^\pm = \Delta x^\pm / \Delta x$, the Shortley-Weller stencil for case (c) is
\begin{equation}\label{eq:shortley_dirichlet_coeffs}
    \partial_x^2 u_i = \frac{1}{\Delta x^2} \left[ \frac{2 u_w^-}{\psi^- (\psi^- + \psi^+)} - \frac{2 u_i}{ \psi^- \psi^+} + \frac{2 u_w^+}{\psi^+ (\psi^- + \psi^+)} \right].
\end{equation}
To recover case (a) let $\psi^+ = 1$ and $u_w^+ = u_{i+1}$, and likewise for case (b) let $\psi^- = 1$ and $u_w^- = u_{i - 1}$ These stencils are well-known and uniquely determined by requiring first order truncation error at all irregular points. The same requirement determines a similar set of coefficients for Neumann boundary conditions, given by
\begin{subequations}\label{eq:shortley_neumann_coeffs}
\begin{align} %
    \partial_x^2 u &= \frac{1}{\Delta x^2} \left[ -\frac{2 \Delta x \partial_x u_w^-}{(2 \psi^- + \psi^+)} - \frac{2 u_i}{\psi^+ (2\psi^- + \psi^+)} + \frac{2 u_{i + 1}}{\psi^+ (2 \psi^- + \psi^+)}\right]. \\
    \partial_x^2 u &= \frac{1}{\Delta x^2} \left[\frac{2 u_{i - 1}}{\psi^+ (\psi^- + 2\psi^+)}- \frac{2 u_i}{ \psi^+ (\psi^- + 2\psi^+)} + \frac{2 \Delta x \partial_x u_w^+}{(\psi^- + 2 \psi^+)} \right]\\
    \partial_x^2 u_i &= \frac{1}{\Delta x} \left[\frac{- \partial_x u_w^- + \partial_x u_w^+ }{\psi^- + \psi^+}\right]
\end{align}
\end{subequations}
For an interface with prescribed jump conditions $[u]$ and $[\beta \partial_x u]$ we take a different approach. Let $x_w \in [x_i, x_{i+1}]$ be the interface location, and let $\Delta x^- = x_w - x_{i}$ and $\Delta x^+ = x_{i+1} - x_w$ be the distance to the nearest grid point on either side. Likewise let $u_w^-$ and $u_w^+$ be the unknown value of $u(x)$ on either side of the interface. The jump conditions are discretized with the first order approximation
\begin{equation} \label{eq:shortley_interface_discretization}
    [u]_w = u_w^- - u_w^+ \quad \text{and} \quad [\beta \partial_x u]_w = \beta^+ \frac{u_{i +1} - u_w^+}{\Delta x^+} - \beta^- \frac{u_w^- - u_i}{\Delta x^-},
\end{equation}
which can be solved to yield the unknown values. Letting $\alpha^\pm = \psi^\pm / \beta^\pm$, 
\begin{equation} \label{eq:shortley_interface_coeffs}
    u_w^+ = \bar{u}_w + \frac{\alpha^+ [u]_w}{\alpha^+ + \alpha^-}, \quad u_w^- = \bar{u}_w - \frac{\alpha^- [u]_w}{\alpha^+ + \alpha^-}, \quad \text{with} \quad \bar{u}_w \equiv \frac{\alpha^+ u_i + \alpha^- u_{i+1}}{\alpha^+ + \alpha^-} - \frac{\Delta x \alpha^+ \alpha^- [\beta \partial_x u]_w}{\alpha^+ + \alpha^-}.
\end{equation}
For homogeneous jump conditions the above simplifies to $u_w^+ = u_w^- = (\alpha^+ u_i + \alpha^- u_{i+1})/(\alpha^+ + \alpha^-)$. These approximate wall values are used in place of a Dirichlet boundary condition when constructing stencils at the irregular points $x_i$ and $x_{i+1}$.

For 2D and 3D discretizations, the irregular stencils of the previous paragraph are applied without alteration along each dimension. For Dirichlet boundary conditions this extension is standard, and the resulting discretization has first order truncation error near domain boundaries. For Neumann conditions or jump conditions we make the approximations $\frac{\partial u}{\partial n} \approx \mathrm{sign}(n_i) \frac{\partial u}{\partial x_i}$ and $\left[ \beta\frac{\partial u}{\partial n}\right] \approx \left[\beta \frac{\partial u}{\partial x_i}\right]$ at points where the boundary or interface intersects a grid line parallel to the $x_i$ axis. The resulting discretizations are not consistent, but as we will demonstrate below they are effective and geometrically robust preconditioners for higher-order immersed discretizations. 

For some immersed geometries, intersection points can be located arbitrarily close to the nearest grid point in the problem domain. In this case the value of $\psi^\pm$ used in Eqs.~\ref{eq:shortley_dirichlet_coeffs},~\ref{eq:shortley_neumann_coeffs},~and~\ref{eq:shortley_interface_coeffs} become arbitrarily small, leading to extremely large coefficients and arbitrarily large eigenvalues in the spectrum of the discretization. As a result, the Shortley-Weller discretizations as outlined above are poor preconditioners for the high order immersed discretizations presented in section~\ref{sec:method}, which have a spectrum that is bounded and scales with $\Delta x^2$. To remedy this issue, the location of each close intersection in the Shortley-Weller discretization is shifted to ensure that $\psi^\pm$ remains bounded away from zero. Specifically, intersections lying on a domain boundary are left in place if the distance between the intersection and the nearest grid point in the problem domain is greater than $\Delta x / 2$. Intersections closer than $\Delta x / 2$ are shifted along the grid line until the distance to their neighboring grid point is exactly $\Delta x / 2$. As a result, $\psi^\pm \ge 0.5$ when computing the coefficients in Eqs.~\ref{eq:shortley_dirichlet_coeffs}~and~\ref{eq:shortley_neumann_coeffs}.
Intersections lying on an immersed interface are shifted to be halfway between the two neighboring grid points, so that $\psi^+ = \psi^- = 0.5$ when computing the coefficients in Eq.~\ref{eq:shortley_interface_coeffs}. This shift reduces the accuracy of the Shortley-Weller discretization, but leads to better preconditioners for our high order immersed discretizations.

\subsection{Shortley-Weller preconditioners for high order methods}
In this section we evaluate the performance of the Shortley-Weller discretization as a preconditioner for higher order immersed methods, and compare its performance with other low-order discretization preconditioners that provide less geometric flexibility. To do so, we choose a 2D test case on a uniform grid, so that the matrix representing each low order discretization can be explicitly assembled, factored, and used as a preconditioner in an iterative linear solver. We report the number of iterations required to reach a solution of a fixed high order discretization as a measure of the effectiveness of each preconditioner. In practice the factorization can be replaced by a less expensive and more scalable multigrid V-Cycle based on the Shortley-Weller discretization, as discussed in the next section.

The geometry and manufactured solutions used in this section are identical to those in section~\ref{subsec:poisson_error}. Given a spatial resolution of $N_x$ points along each axis, we assemble the linear system $\vb{L}_{\Omega} \vb{u}_{\Omega} = \vb{f}_\Omega - \vb{L}_\Gamma \vb{u}_\Gamma$ using a $(4,5)$ high order immersed discretization. The system is solved with left-preconditioned GMRES, where the preconditioner is a sparse LU factorization of a low order immersed discretization. Specifically, we choose as preconditioners the Shortley-Weller discretization as discussed above, two second order stencils with third- and fourth-order IIM interpolants, respectively denoted $(2,3)$ and $(2,4)$, and two fourth order stencils with third and fourth order IIM interpolants, respectively denoted $(4,3)$ and $(4,4)$. In each case, we consider the solution converged when the $L_2$ norm of the preconditioned residual reaches a relative tolerance of $\epsilon_r = 10^{-6}/N_x$ compared to the $L_2$ norm of the right hand side. 

\begin{figure}[tb!]
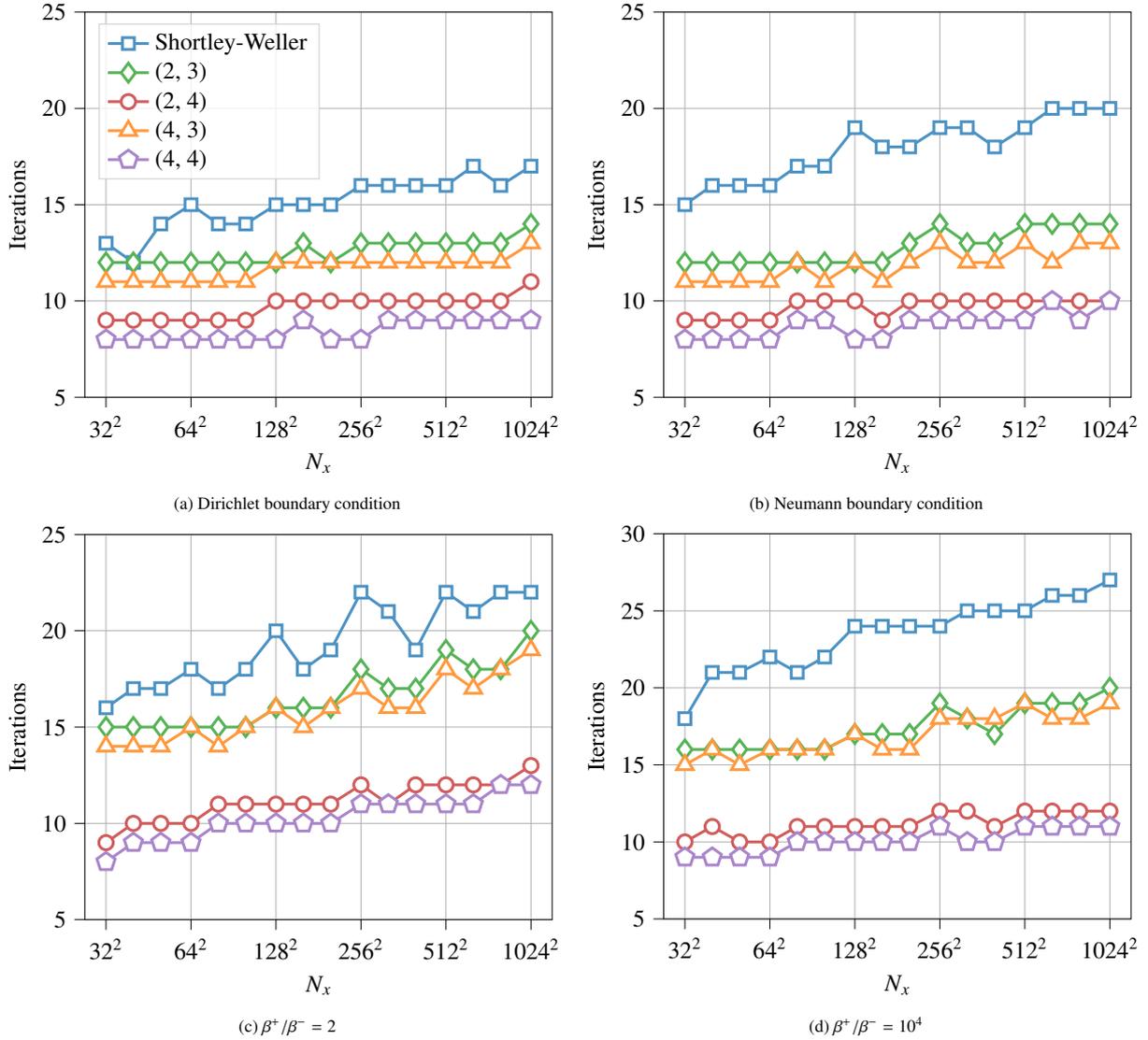

    \centering
    \begin{subfigure}{0.49\textwidth}
        \centering
        \resizebox{\textwidth}{!}{\input{results/stencil/low_order_prec_dirichlet_form}}
        \caption{Dirichlet boundary condition}
        \label{subfig:lower_order_dirichlet}
    \end{subfigure}
    \begin{subfigure}{0.49\textwidth}
        \centering
        \resizebox{\textwidth}{!}{\input{results/stencil/low_order_prec_neumann_form}}
        \caption{Neumann boundary condition}
        \label{subfig:lower_order_neumann}
    \end{subfigure}
        \centering
    \begin{subfigure}{0.49\textwidth}
        \centering
        \resizebox{\textwidth}{!}{\input{results/stencil/low_order_prec_flux_low_form}}
        \caption{$\beta^+ / \beta^- = 2$}
        \label{subfig:lower_order_low_ratio}
    \end{subfigure}
    \begin{subfigure}{0.49\textwidth}
        \centering
        \resizebox{\textwidth}{!}{\input{results/stencil/low_order_prec_flux_high_form}}
        \caption{$\beta^+ / \beta^- = 10^4$}
        \label{subfig:lower_order_high_ratio}
    \end{subfigure}
    \caption{Iterations required to solve the $(4,5)$ discretization of the Poisson equation using various low order preconditioners with (a) Dirichlet boundaries, (b) Neumann boundaries, (c) a material interface with $\beta^-/\beta^+ = 2$, and (d) a material interface with $\beta^-/\beta^+ = 10^4$.
    }
    \label{fig:lower_order_iterations}
\end{figure}

Figures~\ref{subfig:lower_order_dirichlet}~and~\ref{subfig:lower_order_neumann} report the iterations required to solve the Poisson equation with Dirichlet or Neumann boundary conditions over a range of spatial resolutions. In each case the number of iterations grows only slightly as the spatial resolution increases, indicating that each preconditioner successfully circumvents the $\mathcal{O}(\Delta x^{-2})$ condition number of the high order system. Figures~\ref{subfig:lower_order_low_ratio}~and~\ref{subfig:lower_order_high_ratio} report similar results for the material interface test case given in Eq.~\ref{eq:manufactured_interface}, using coefficients $\beta^- = 1$ and $\beta^+ = 1/2$ or $\beta^+ = 10^{-4}$ respectively. The performance of each lower order $(n,k)$ preconditioner decreases slightly as both the interior order $n$ and boundary order $k$ decreases, and the Shortley-Weller preconditioner is the least performant across all boundary and interface conditions. However, the performance of the Shortley-Weller preconditioner is still comparable to the fully consistent $(2,3)$ and $(2,4)$ preconditioners, despite the shifting of boundary points and inconsistent treatment of Neumann and interface conditions. The slight decrease in performance is justified by the greatly increased geometric flexibility of the discretization, which can be applied without issue to the coarse grids of a multigrid algorithm as described in the next section.

\subsection{Multigrid preconditioners based on the Shortley-Weller discretization}

In addition to its geometric robustness and relative simplicity, the Shortley-Weller discretization can also be solved effectively with geometric multigrid. A simple and robust method for homogeneous Dirichlet boundary conditions is given by \citet{Trottenberg2000}. The method consists of coarsening the grid by a factor of two between multigrid levels, applying a red-black Gauss-Seidel smoother on each level, and transferring data between grids with linear interpolation (coarse to fine) or half weighted restriction (fine to coarse). These transfer operations are applied to the full computational grid without accounting for immersed interfaces, with zero values assigned to points outside of the problem domain. Results from \cite{Trottenberg2000} demonstrate rapid convergence: fine-level residuals are reduced by a factor of $\rho_i = \Vert r^{i}\Vert /\Vert r^{i-1} \Vert < 0.09$ for each multigrid iteration $i$ across a range of irregular geometries. We adopt this multigrid approach directly for our Shortley-Weller discretizations in case for Dirichlet boundary conditions.

For the homogeneous Dirichlet case, the data transferred between grids is expected to approach zero near domain boundaries, and assigning zero values to points outside of the problem domain introduces only small errors during interpolation or restriction. For Neumann conditions this is not the case, and satisfactory convergence rates require interpolation and restriction operators that account for immersed boundaries and interfaces. Numerical experiments demonstrate that good performance can be recovered by leaving the half weighted restriction unaltered, and replacing the standard bilinear or trilinear interpolation stencil with an average over all neighboring coarse values that fall within the problem domain. These formulations are equivalent away from the interface, but for points adjacent to the interface the averaging formulation retains $\mathcal{O}(\Delta x)$ accuracy while the linear interpolation leads to $\mathcal{O}(1)$ interpolation errors. For the rare case of fine level points with no valid coarse neighbors, we assign an interpolated value of zero. In the context of a multigrid method, the interpolated field is a correction to the fine level solution, and this treatment delays corrections at these pathological points until the algorithm reaches a finer resolution level that can resolve the local geometry.

For problems with a piecewise constant coefficient $\beta$ and prescribed homogeneous jump conditions, the solution and its derivative are continuous across each immersed interface. This allows standard interpolation and restriction operators to be applied to the solution without introducing large errors near the interface. However applying the multigrid method outlined above to the operator $\nabla \cdot (\beta \nabla u)$ with large jumps in $\beta$ leads to poor accuracy. The issue arises from an imbalance in the residual: when $\beta^+ \gg \beta^-$, even large errors in the solution on $\Omega^-$ can have a negligible impact on residual norm $\Vert r \Vert = \Vert f - \nabla \cdot (\beta \nabla u) \Vert$. To remedy this, the multigrid method is applied to the equivalent constant coefficient problem $\nabla^2 u = f/\beta$ on each domain, using the same interface treatment described in the previous section to determine Dirichlet boundary values on each side of the interface.

\begin{figure}[tb!]
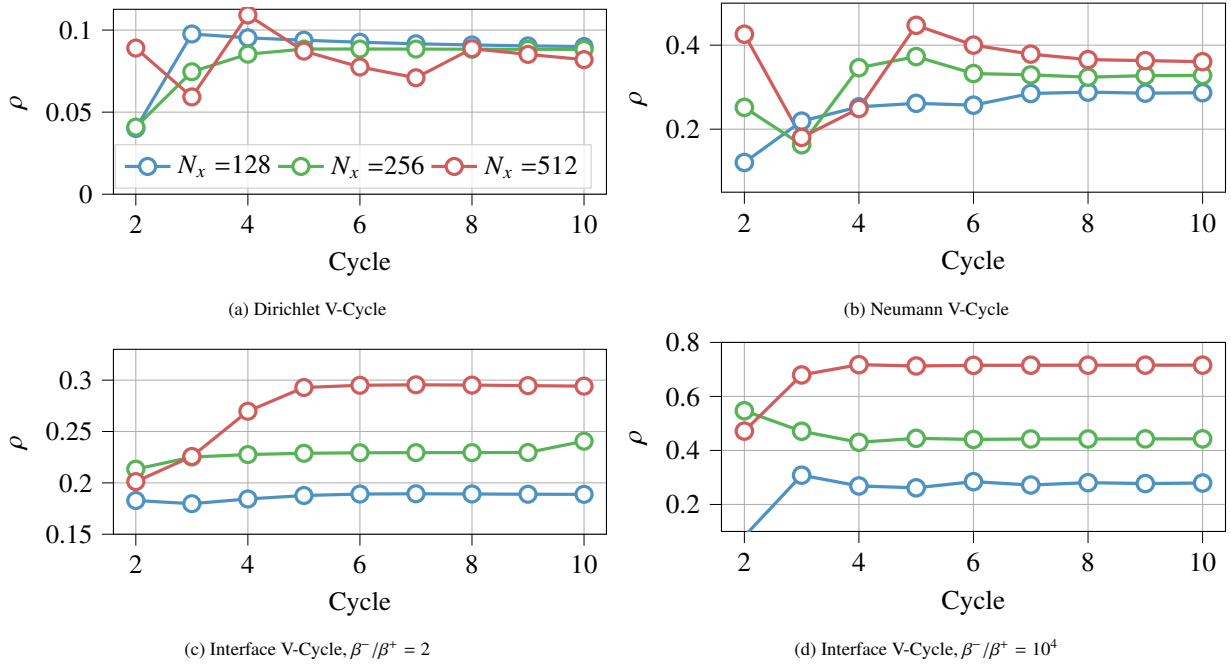

    \centering
    \begin{subfigure}{0.49\textwidth}
        \centering
        \resizebox{\textwidth}{!}{\input{results/stencil/mg_2d_v_convergence_interface_dirichlet}}
        \caption{Dirichlet V-Cycle}
        \label{subfig:mg_2d_dir_v}
    \end{subfigure}
    \begin{subfigure}{0.49\textwidth}
        \centering
        \resizebox{\textwidth}{!}{\input{results/stencil/mg_2d_v_convergence_interface_neumann}}
        \caption{Neumann V-Cycle}
        \label{subfig:mg_2d_neu_v}
    \end{subfigure}
    \begin{subfigure}{0.49\textwidth}
        \centering
        \resizebox{\textwidth}{!}{\input{results/stencil/mg_2d_v_convergence_interface_low}}
        \caption{Interface V-Cycle, $\beta^-/\beta^+ = 2$}
        \label{subfig:mg_2d_low_v}
    \end{subfigure}
    \begin{subfigure}{0.49\textwidth}
        \centering
        \resizebox{\textwidth}{!}{\input{results/stencil/mg_2d_v_convergence_interface_high}}
        \caption{Interface V-Cycle, $\beta^-/\beta^+ = 10^4$}
        \label{subfig:mg_2d_high_v}
    \end{subfigure}
    \caption{Convergence history of multigrid V-Cycles applied to the Shortley-Weller discretization with (a) Dirichlet boundaries, (b) Neumann boundaries, (c) a material interface with $\beta^-/\beta^+ = 2$, and (d) a material interface with $\beta^-/\beta^+ = 10^4$. Results are given for three different spatial resolutions.}
    \label{fig:2d_mg_v_convergence}
\end{figure}

To assess how well each multigrid method approximates the inverse of the corresponding Shortley-Weller discretization, we report the observed convergence rate $\rho_i$ for each multigrid iteration $i$ when applied to the test cases used in the previous section. The multigrid algorithm used is as described above, with a direct coarse solver applied on the level with $N_x = 8$. For each test case the geometry and forcing terms are unchanged, but the boundary conditions or jump conditions are made homogeneous to mimic the action of a preconditioner. We show results for spatial resolutions $N_x = [128,\,256,\,512]$, and note that for an ideal multigrid algorithm the convergence rate is largely independent of spatial resolution. Figure~\ref{fig:2d_mg_v_convergence} plots the observed convergence rates of a multigrid V-Cycle algorithm.
Consistent with \cite{Trottenberg2000}, the Dirichlet case exhibits rapid convergence independent of resolution. For Neumann boundaries or jump conditions the algorithms are convergent, but their performance decreases for larger problem sizes.

\begin{figure}[tb!]
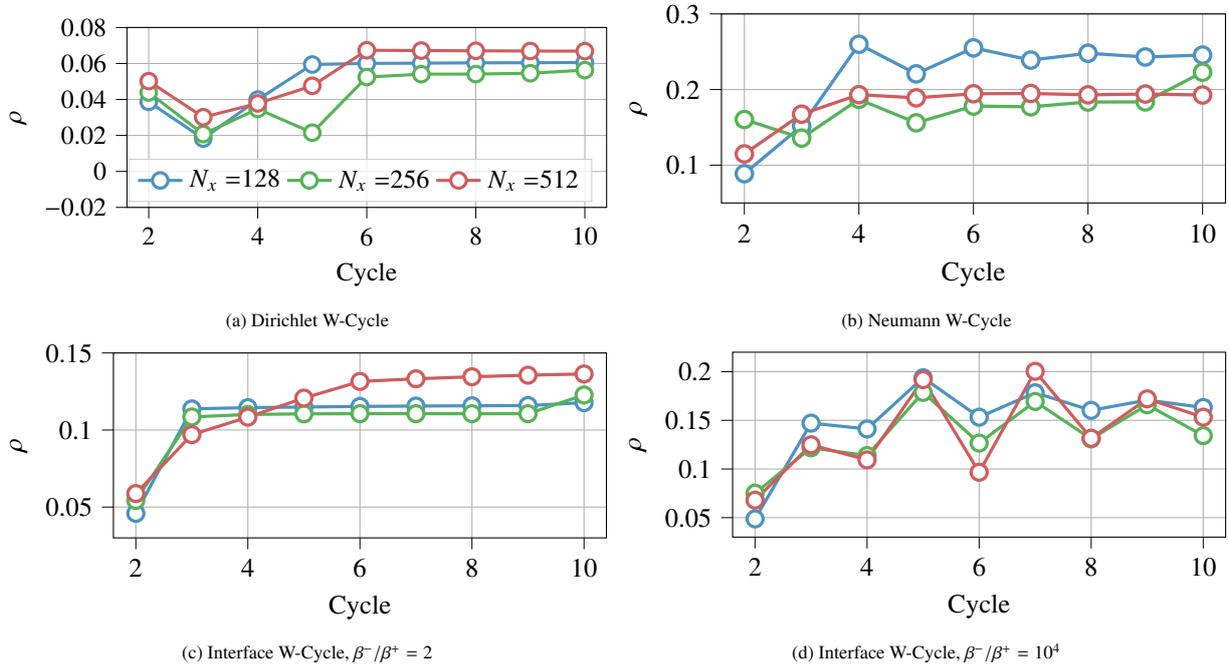

    \centering
    \begin{subfigure}{0.49\textwidth}
        \centering
        \resizebox{\textwidth}{!}{\input{results/stencil/mg_2d_w_convergence_interface_dirichlet}}
        \caption{Dirichlet W-Cycle}
        \label{subfig:mg_2d_dir_w}
    \end{subfigure}
    \begin{subfigure}{0.49\textwidth}
        \centering
        \resizebox{\textwidth}{!}{\input{results/stencil/mg_2d_w_convergence_interface_neumann}}
        \caption{Neumann W-Cycle}
        \label{subfig:mg_2d_neu_w}
    \end{subfigure}
    \begin{subfigure}{0.49\textwidth}
        \centering
        \resizebox{\textwidth}{!}{\input{results/stencil/mg_2d_w_convergence_interface_low}}
        \caption{Interface W-Cycle, $\beta^-/\beta^+ = 2$}
        \label{subfig:mg_2d_low_w}
    \end{subfigure}
    \begin{subfigure}{0.49\textwidth}
        \centering
        \resizebox{\textwidth}{!}{\input{results/stencil/mg_2d_w_convergence_interface_high}}
        \caption{Interface W-Cycle, $\beta^-/\beta^+ = 10^4$}
        \label{subfig:mg_2d_high_w}
    \end{subfigure}
    \caption{Convergence history of multigrid W-Cycles applied to the Shortley-Weller discretization with (a) Dirichlet boundaries, (b) Neumann boundaries, (c) a material interface with $\beta^-/\beta^+ = 2$, and (d) a material interface with $\beta^-/\beta^+ = 10^4$. Results are given for three different spatial resolutions. Unlike the V-Cycle, the spatial resolution has only a small impact on the convergence factor.}
    \label{fig:2d_mg_w_convergence}
\end{figure}

Figure~\ref{fig:2d_mg_w_convergence} reports the result of applying a multigrid W-cycle to the same test cases, keeping the same smoothing, restriction, and interpolation operators unchanged. The W-cycle leads to faster convergence for each spatial resolution and boundary condition, and unlike the V-cycle the convergence rate is largely independent of spatial resolution. The Dirichlet multigrid algorithm still shows the fastest convergence, with the Neumann algorithm converging slowest and the two immersed interface test cases falling somewhere in between.

\begin{figure}[tb!]
    \centering
        \resizebox{0.49\textwidth}{!}{\input{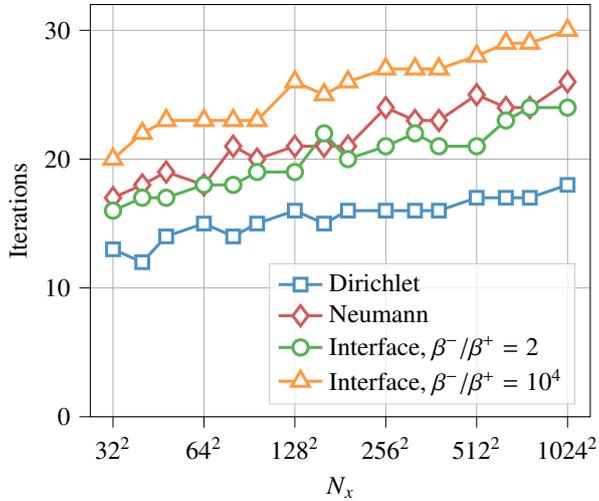}}
    \caption{Iterations required to solve the $(4,5)$ discretization of the Poisson equation using multigrid preconditioning. Results include Dirichlet and Neumann boundary conditions, as well as prescribed jump conditions.}
    \label{fig:2d_mg_iterations}
\end{figure}

As a final test, we solve the same Poisson problems used for the results in Fig.~\ref{fig:lower_order_iterations}, replacing the exact inverse of the Shortley-Weller discretization with a single multigrid V-Cycle. Figure~\ref{fig:2d_mg_iterations} plots the number of iterations required for each boundary condition across a range of spatial resolutions. As with the exact inverse, the number of iterations varies only weakly with spatial resolution, indicating that the multigrid preconditioner effectively mitigates the $\mathcal{O}(\Delta x^{-2})$ condition number of the high order immersed discretization. While the W-cycle leads to faster convergence for the multigrid methods shown in Fig.~\ref{fig:2d_mg_w_convergence}, we observed only a slight improvement the performance of the preconditioner when replacing the single V-Cycle with a single W-cycle. For efficiency the single multigrid V-Cycle is used as a preconditioner in the remainder of this work.

\section{Results and performance on 3D problems}\label{sec:results}

In this section we will demonstrate the presented approach in the context of 3D problems. First, we will describe the 3D implementation in the multiresolution adaptive grid solver \textsc{murphy} \cite{Gillis2022}. Second, we will demonstrate convergence of the solution, and analyze the convergence of the iterative solver in the 3D implementation. Third, we will show the parallel scalability of the resulting algorithm. Finally, we consider a set of applications in 3D relying on a combination of complex 3D geometries and multiresolution grid adaptation to solve practical problems. All results in this section will be limited to immersed boundary problems, following the same algorithm as described above.  

\subsection{Implementation in 3D}
The 3D implementation relies on \textsc{murphy}, a block-structured octree grid framework relying on interpolating wavelets to perform adaptive grid simulations of partial differential equations \cite{Gillis2022}. In turn, \textsc{murphy} builds upon the octree capabilities of \textsc{p4est} \cite{Burstedde:2011}. 

For solving the Poisson equation in 3D domains, we interfaced \textsc{murphy} with the matrix-free algorithms implemented in the PETSc library \cite{petsc-web-page}. In particular, the matrix-vector product of the high order discretization stencil is implemented as a Laplacian stencil operator with immersed interface corrections, identical to the 2D approach presented above \cite{Gabbard2024}. For multiresolution grids, wavelet transforms are used to create same-level ghost layers across resolution boundaries, after which standard Laplacian stencils can be used in the interior of each block \cite{Gillis2022}. We enforce that all blocks containing an immersed interface control point, and all of their neighbors, are refined to a user-specified maximum resolution. This allows all IIM operations to be done at locally uniform resolution \cite{Gabbard2024}. The other grid blocks are refined or coarsened based on the wavelet transform of the right-hand side of the Poisson equation, using wavelet orders  consistent with the discretization order of the finite difference schemes \cite{Gillis2022}.

The 3D IIM multigrid preconditioner is based on a free space Full Approximation Storage multigrid solver implemented in \textsc{murphy} \cite{Poncelet2025}, following \cite{Trottenberg2000}. We implement custom prolongation and restriction operators, as well as the low order Shortley Weller stencils, consistent with the 2D implementation described above. Since the smoothing in \textsc{murphy} occurs necessarily in a block-based manner, the red-black Gauss-Seidel scheme of the 2D uniform grid implementation is replaced by a block-based Gauss-Seidel scheme with lexicographic ordering within each block. Inside the multigrid hierarchy in \textsc{murphy}, the grid is coarsened to a single block in each unit cube domain. At this coarsest level the grid resolution is uniform, and we can form the matrix associated with the Shortley-Weller stencil of the preconditioner. This matrix is passed to the coarse level solver, which in our case is a GMRES solver left-preconditioned with the Geometric Algebraic Multigrid (GAMG) \cite{petsc-user-ref}. Because this coarse grid solver in the preconditioner is iterative, it violates the linearity of the preconditioner \cite{Saad:2003}. Consequently, we found that using a right-preconditioned Flexible GMRES as the high order fast iterative solver \cite{Saad:1993}, as opposed to the standard GMRES used in the 2D context, significantly improves the robustness of our 3D algorithm. To alleviate high memory overheads, we further use a FGMRES restart frequency of 10 throughout the tests below.

\subsection{Verification and analysis for immersed boundaries}
\label{subsec:3D_boundaries}
\begin{figure}[tb!]
    \centering
    \begin{subfigure}{0.43\textwidth}
        \centering
\includegraphics[width=\textwidth]{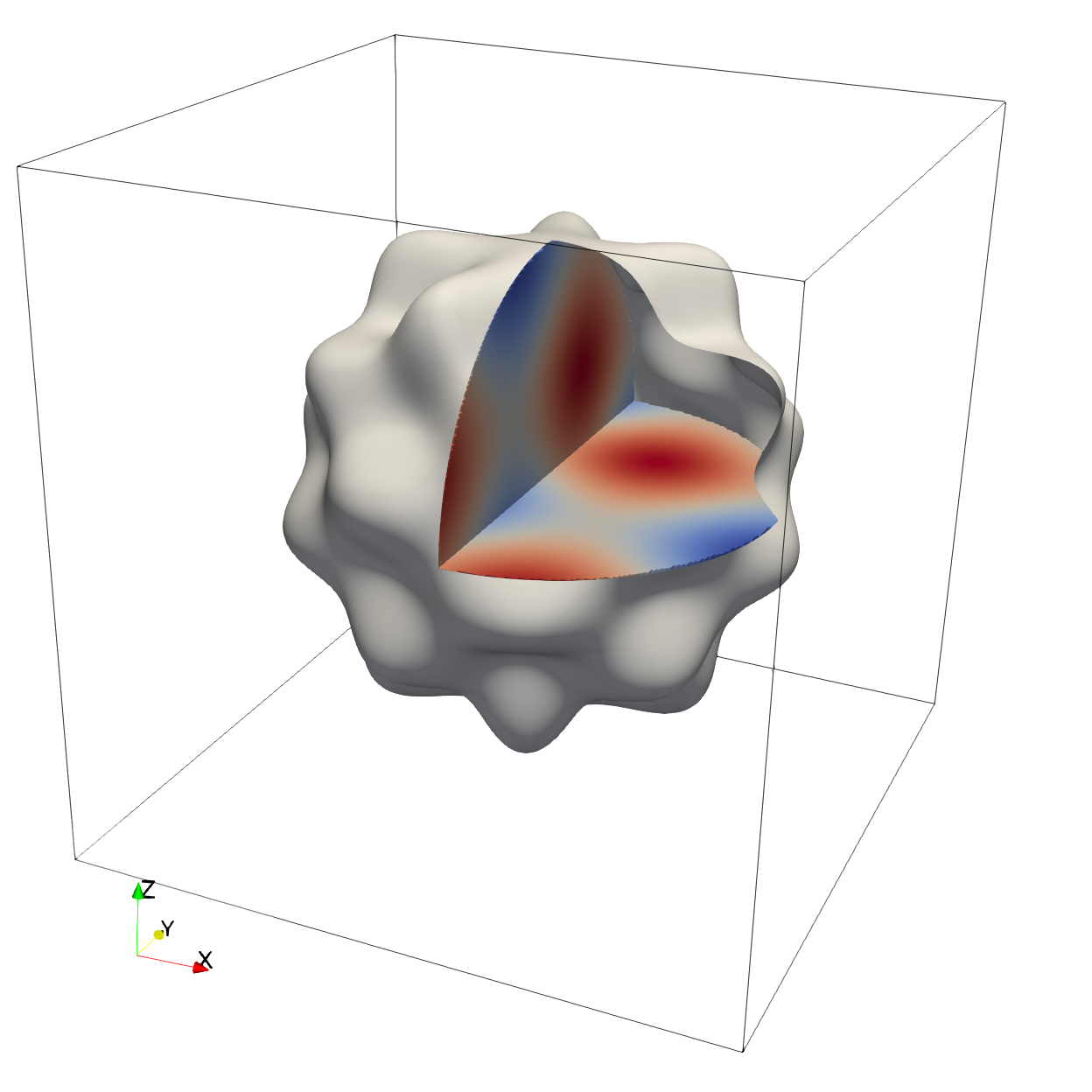}
        \caption{Geometry and exact solution}
        \label{subfig:conv3d_setup}
    \end{subfigure}
    \hfill
    \begin{subfigure}{0.53\textwidth}
        \centering
        \resizebox{\textwidth}{!}{\input{results/stencil/3D_convergence_all}}
        \caption{Convergence of maximum solution error across different discretization orders}
        \label{subfig:conv3d_error}
    \end{subfigure}
    \caption{Convergence of an exterior, 3D Poisson equation with a manufactured solution inside a perturbed sphere. We vary the $(n,k)$ discrestization order, where $n$ is the order of the finite-difference scheme and $k$ the order of the IIM polynomial interpolant.}
    \label{fig:conv3d}
\end{figure}

To test convergence and performance of  problems with immersed boundaries, we set up a manufactured problem whose solution is
\[
u(\xv) = \sin(4 \pi x_1) + \sin(4 \pi x_2) + \sin(4 \pi x_3).
\]
The Poisson problem is stated within the interior of a sphere that is radially perturbed by a distance proportional to a spherical harmonic, as in \cite{Gabbard2024}. The sphere has radius $\rho_0 = 0.35$, perturbation amplitude $\tilde{\rho} = 0.05$, and is embedded within a unit cube domain (Fig.~\ref{subfig:conv3d_setup}). The right-hand side and boundary conditions (Dirichlet or Neumann) are evaluated analytically from the exact solution and used as input to the solver. To verify convergence, we set the relative tolerance of the iterative solver to $10^{-15}$ and perform a convergence test on uniform grids with resolutions $48^3$ to $768^3$. At each resolution, we evaluate the maximum error between the exact and numerical solution over all the grid points within the perturbed sphere.  Figure~\ref{subfig:conv3d_error} shows the convergence of the maximum solution error computed using Dirichlet and Neumann boundary conditions. The results demonstrate that the proposed high order immersed approach can be translated directly to large-scale problems in 3D.

\begin{figure}[htb!]
    \centering
    \resizebox{0.49\textwidth}{!}{\input{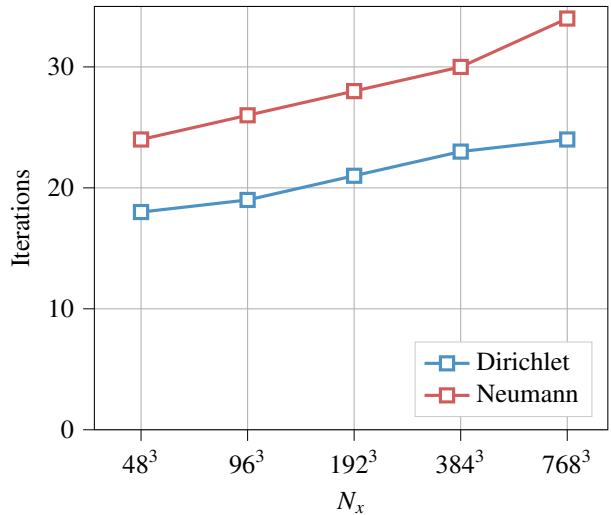}}
    \label{subfig:errorConv_sc}
     \caption{Iterations required to solve the $(4,5)$ discretization of the Poisson equation in 3D, using the proposed Shortley-Weller multigrid preconditioning. Results include both Dirichlet and Neumann boundary conditions.}
    \label{fig:3D_iterative_convergence}
\end{figure}

To assess the effectiveness of the proposed SW-MG preconditioning strategy, we record the number of iterations required to reduce the relative residual to $10^{-6}/N_x$ as a function of spatial discretization order, resolution, and type of boundary condition, similar to Fig~\ref{fig:2d_mg_iterations} in 2D. The results in 3D, shown in Figure~\ref{fig:3D_iterative_convergence}, are similar to those in 2D: the algorithm is robust and the residual converges consistently and at roughly similar rates across different resolutions. The highest iteration count is recorded for the high order Neumann  demonstrate that for either order of accuracy, the Neumann condition requires more iterations. This can be attributed to the fact that the Neumann Shortley-Weller preconditioner is not a consistent discretization to the low order problem. Nevertheless, we find that the algorithm is robust and the residual converges consistently and at roughly similar rates across different resolutions.

\subsection{Parallel performance}
To test the parallel performance of the 3D implementation, we perform a parallel scaling test on the Cray HPE EX supercomputer Perlmutter, part of the National Energy Research Scientific Computing Center (NERSC). We utilize Perlmutter's CPU-only nodes, each consisting of two AMD EPYC 7763 CPUs for a total of 128 physical compute cores and 512 GB DDR4 memory per node. The interconnect is HPE Slingshot 11. On this architecture, \textsc{murphy} is compiled using the Cray MPICH compiler suite version 8.1.25, based on the GNU C++ compilers version 12.2.0. 

We setup the parallel scaling study to solve the Poisson equation in a periodic domain, exterior to an immersed interface. The interface and the problem set are identical to those described in section~\ref{subsec:3D_boundaries} and shown in Fig~\ref{subfig:conv3d_setup}. We choose a uniform grid resolution, fourth-order discretization, and Dirichlet  boundary conditions for all problems. The Poisson equation is solved to relative residual of \num{1e-10}, requiring 29 iterations and yielding an $L_\infty$ error of $\num{7e-9}$ for all problem sizes. 

\paragraph{Weak scaling} Starting with weak scaling, the reference domain has size $[1, 1, 2]$ with resolution $2\cdot768^3 \approx \num{0.9e9}$ grid points, distributed across two nodes. When increasing the node count, periodic images are added in the $z$ direction so that the domain size becomes $(1,1,N)$ for $N$ nodes. Note that when adding periodic images of the domain, we also include the immersed boundary in the image. We report scaling results up to \num{256} nodes (\num{32768} cores), with a total number of about 116 billion unknowns. The results (Fig.~\ref{fig:weak_scaling}) demonstrate that the total time-to-solution only mildly increases with node count. With the weak scaling efficiency defined as $\eta_w(N) = T(2)/T(N)$, where $T(i)$ is the time-to-solution on $i$ nodes, we find that the solver maintains 95\% overall efficiency even at $N=256$. From the timing plot we observe that the multigrid preconditioner is responsible for about half the total cost, following by the matrix-vector multiplication (\textit{MatVec}), and then the FGMRES algorithm. All three components scale similarly.

\paragraph{Strong scaling}
For the strong scaling study, we keep the domain size fixed at $2\cdot 768^3$ at 2 nodes, and subsequently increase the node count in factors of 2 up to $N=64$. On 2 nodes, we then process about \num{3.5e6} unknowns per core, whereas on 64 nodes this number drops to only about \num{55e3}. Figure~\ref{fig:strong_scaling} shows the result in terms of time-per-iteration in panel (a) as well as strong scaling efficiency defined as $\eta_s(N) = (2/N) T(2)/T(N)$ in panel (b). The results show that the matrix-vector multiplication as well as the FGMRES algorithm scale perfectly, whereas the multigrid preconditioner scales the worst. Nevertheless, we maintain an 80\% strong efficiency when increasing from 2 nodes to 16 nodes, associated with a speedup of $6.4 \times$. After this the scaling performance deteriorates, with a speed up of only about $16\times$ at 64 nodes compared to 2 nodes, so that the strong efficiency drops to about 50\%. Overall the deterioration of strong scaling as the number of unknowns per core decreases is common to elliptic solvers, as communication overheads start to dominate. In our specific algorithm, furthermore, it is expected that the multigrid solver poses the bottleneck in strong scaling, as the inherent multiscale character of the algorithm is associated with significant load imbalances at the coarser levels. This is acceptable at high numbers of unknowns per core, as most time is spent in the finest level smoothing and prolongation/restriction steps. However, as the number of unknowns per core decreases, the load imbalances on the coarser levels become more and more dominant, which deteriorates the strong scaling. Future work will be focused on improving the performance, e.g.\ by exploring a reduction in the depth of the multiscale hierarchy when the number of unknowns per core is small. 

\begin{figure}[htb!]
    \centering
        \begin{subfigure}{0.49\textwidth}
        \centering
    \resizebox{\textwidth}{!}{\input{results/stencil/weak_scaling}}
    \caption{Timings}
    \label{subfig:weak_scaling_time}
    \end{subfigure}
     \begin{subfigure}{0.49\textwidth}
         \centering
         \resizebox{\textwidth}{!}{\input{results/stencil/weak_efficiency}}
         \caption{Efficiency}         \label{subfig:weak_scaling_efficiency}
     \end{subfigure}
     \caption{Weak scaling results from 2 nodes (128 cores, \num{0.9e9} grid points) to 256 nodes (\num{32768} cores, \num{116e9} grid points): time-per-iteration (a) and weak scaling efficiency $\eta_w$ (b) against the number of nodes. The sum of the  timings (black) are split between the multigrid preconditioner (red), the high order matrix-vector multiplication (green), and the FGMRES algorithm (blue).}
    \label{fig:weak_scaling}
\end{figure}

\begin{figure}[htb!]
    \centering
        \begin{subfigure}{0.49\textwidth}
        \centering
    \resizebox{\textwidth}{!}{\input{results/stencil/strong_scaling}}
    \caption{Timings}
    \label{subfig:strong_scaling_time}
    \end{subfigure}
     \begin{subfigure}{0.49\textwidth}
         \centering
         \resizebox{\textwidth}{!}{\input{results/stencil/strong_efficiency}}
         \caption{Effiency}         \label{subfig:strong_scaling_efficiency}
     \end{subfigure}
     \caption{Strong scaling results from 2 nodes (128 cores, \num{3.5e6} grid points per core) to 64 nodes (\num{8192} cores, \num{55e3} grid points per core): time-per-iteration (a) and weak scaling efficiency $\eta_s$ (b) against the number of nodes. The sum of the  timings (black) are split between the multigrid preconditioner (red), the high order matrix-vector multiplication (green), and the FGMRES algorithm (blue).}
    \label{fig:strong_scaling}
\end{figure}

\subsection{Applications in 3D}
To demonstrate the ability of the algorithm and implementation to solve elliptic problems in complex 3D domains, we consider two applications. First, we demonstrate a heat transfer application inside a triply periodic gyroid structure. Subsequently, we consider an exterior potential flow problem around a set of interlinked tori forming a larger scale ring structure. Both cases are solved using the high order discretization present above, with the latter also exploiting the multiresolution adaptive grid.

\subsubsection{Heat conductivity inside a triply periodic structure}
Here we consider a Poisson problem representative of a heat conductivity test inside a triply periodic gyroid domain. The domain is prescribed by the level set function
\[
\phi(\vb{x}) = \cos(k_1 x_1) \sin(k_2 x_2) + \cos(k_2 x_2) \sin(k_3 x_3) + \cos(k_3 x_3) \sin(k_1 x_1) - \alpha, 
\]
where here we set $(k_1,k_2,k_3) = (8\pi,8\pi,8\pi)$ and $\alpha = -0.5$. The domain is defined as the region $\phi(\vb{x})> 0$ bounded by the level-set $\phi(\vb{x}) = 0$, solved inside a unit cube with periodic boundary conditions in all directions.

\begin{figure}[htb!]
    \centering
    \includegraphics[width=\textwidth]{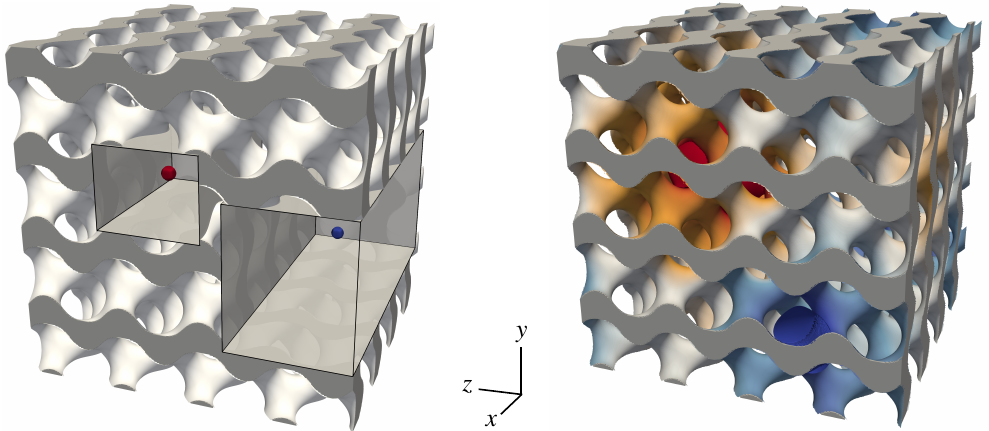}
    \caption{Right-hand side for the heat equation (left) consisting of two isolated sources of opposite sign in different locations of the periodic gyroid structures. Parts of the structure are cut out to better show the location of the sources. Solution for the temperature (right) shows the temperature field associated with these sources. In both plots, the heat equation is solved in the open parts of the visualization; the `outside' of the gyroid structure is masked with dark grey.}
    \label{fig:gyroid}
\end{figure}

We solve the non-dimensional heat equation $\nabla^2 u = -q$ where $u$ is the non-dimensional temperature, $q$ the non-dimensional heat flux, and we assume a  constant conductivity. Here $q$ consists of two opposite signed sources embedded inside the gyroid structure illustrated in Fig~\ref{fig:gyroid}. Both sources are described by compact isotropic Gaussian profiles 
\[
f(r) = \begin{cases}
    \exp\left(\frac{-(r/\sigma)^2}{1 - (r/\gamma)^2}\right), & (r/\gamma)^2 < 1 \\
    0, & \text{else}
\end{cases}
\]
where here the support $\gamma = 1/40$ and $\sigma = \gamma/2$. The source field is then prescribed as $q(\vb{x}) = \alpha^+ f(|\vb{x} - \vb{x}_q^+|) + \alpha^- f(|\vb{x} - \vb{x}_q^-|)$ with $\alpha^\pm = \pm 10^5$ for the positive and negative source, respectively. The sources are located at $\vb{x}_q^- = (3/32, 7/32, 9/32)$ and $\vb{x}_q^+ = (23/32, 19/32, 21/32)$. On the immersed boundary we impose homogeneous flux boundary conditions. The problem is solved using a fourth order discretization on a uniform domain grid with resolution $384^3$ to a relative residual of \num{1e-10}. 

Figure~\ref{fig:gyroid} shows the right-hand side field (left) and the temperature solution (right) inside the lattice. The part of the gyroid that is not solved for (i.e.\ where $\phi(\vb{x}) < 0$) is masked as a dark grey solid. For the solution (right), we visualize the temperature field on the boundary of the gyroid. The solution shows the spread of temperature due to the isolated sources and their periodic neighbors. The example highlights the solver's ability to solve elliptic problems in complex geometries, and recover to high order both the solution as well as surface quantities -- in this case the surface temperature field. 

\subsubsection{Potential flow around a ring of interlinked tori}
We evaluate the potential flow problem external to an immersed geometry. The geometry is constructed of a set of 24 connected tori forming a ring structure. The rings are immersed in a unit cube domain. Within the domain, the radius of the ring is $0.3$, Each torus has major radius \num{5.57e-2} and minor radius \num{1.19e-2}, so that the edge-to-edge distances are maintained at exactly $0.01$. We solve the Laplace problem for the potential function $\Phi$, with non-homogeneous Neumann boundary condition on the immersed boundary $\partial_n \Phi = -\vb{U_\infty} \cdot \vb{n}$, where here $\vb{U}_{\infty} = (1, 0, 0)^T$. On the domain boundaries we assign  homogeneous Neumann conditions. The velocity field is then recovered as $\vb{u} = \nabla \Phi + \vb{U}_\infty$, so that on the immersed surface $\vb{u} \cdot \vb{n} = 0$. Both the Laplace equation as well as the subsequent gradient are discretized using fourth order central finite difference schemes, yielding an overall fourth order solution for the velocity field as demonstrated in section~\ref{subsec:poisson_error}.

\begin{figure}[htb!]
    \centering
    \includegraphics[width=\textwidth]{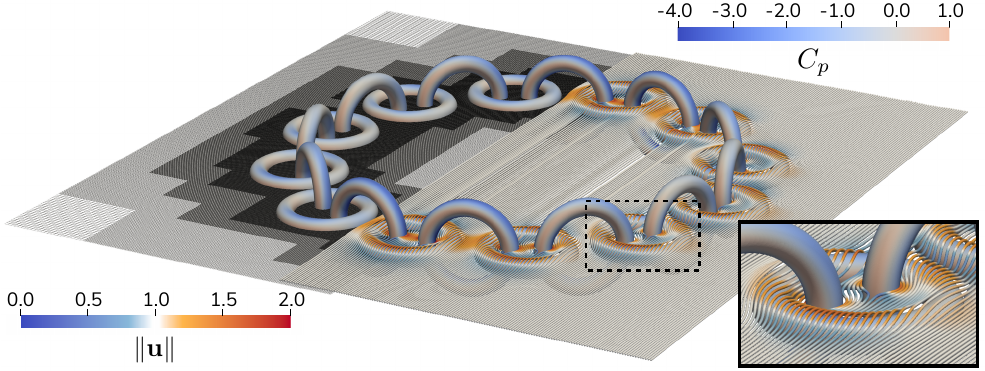}
    \caption{Potential flow solution around a ring of interlinked tori. The streamlines are colored by velocity magnitude $\| \vb{u}\|$, and the rings are colored by pressure coefficient $C_p$. The flow direction is from bottom left to top right, or vice versa. The left half of the plot shows the multiresolution adapted grid, and the inset on the bottom-right shows a close up of the flow around one of the rings (dashed rectangle in the main figure).}
    \label{fig:torusring}
\end{figure}

The problem is solved on a multiresolution adapted grid at maximum level 6 ($262144$ blocks) with a linear blocksize of $N=24$ grid points. The effective resolution (if the grid was refined to the maximum level everywhere) is thus $1536^3 \approx \num{3.6e9}$ grid points. The grid adaptation policy is to refine each grid block containing an immersed boundary control points, as well as all of its neighbors, to the finest level. Since the right-hand side in this problem is zero, the rest of the grid is coarsened as much as possible, while maintaining the overall 1:2 resolution jump condition described in \cite{Gillis2022}. After performing the grid adaptation, the actual grid contains $16206$ blocks for a total number of grid points of $\num{0.22e9}$, a reduction by a factor of about 16. We use 360 compute cores to solve  the Laplace equation for the potential field to a relative residual of $\num{1e-10}$. Figure~\ref{fig:torusring} shows a visualization of the streamlines of $\vb{u}$, as well as a cross-section of the adapted grid. We further compute the velocity vector on each control point by evaluating the wall derivatives of $\Phi$ on the interface. From this we compute the surface pressure coefficient as $C_p = 1 - \|\vb{u}\|^2/\|\vb{U}_\infty\|^2$, which is visualized on the surface of the rings in Figure~\ref{fig:torusring}.

\section{Conclusion}\label{sec:conclusion}
Our work provides the foundation for the use of high order immersed methods for elliptic problems in large scale 3D simulations. The proposed low-order Shortley-Weller stencil in the multigrid preconditioner is robust to arbitrary immersed geometries, allowing the multigrid to operate throughout the entire grid hierarchy. We demonstrate that the multigrid acts as an effective preconditioner for a GMRES solver of the full high order immersed discretization, so that the number of iterations does not significantly depend on grid resolution, boundary conditions, or coefficient jumps across the immersed interface. As the entire algorithm can be implemented in a matrix-free manner, the extension to a large-scale parallel multiresolution 3D solver is relatively straightforward. The resulting approach provides a unified way to tackle 3D Poisson problems with high order discretizations of boundary and/or interface conditions imposed on immersed geometries.

The algorithm can be further improved in various ways, which we reserve for future work. First, though the treatment of non-constant coefficients in the high order stencil is straightforward, integrating this into the multigrid solver requires further work. Second, further work could potentially reduce the iteration count by (1) improved preconditioning of the augmented system for singular problems \cite{Min2006}; (2) investigation into multigrid prolongation/restriction operators adapted more specifically to the Shortley-Weller stencil; (3) further hyperparameter and algorithm selection tuning for the multigrid and the GMRES solver. Specifically for the 3D implementation, we can further improve the computational cost by reducing some spurious ghosting calls associated with multiresolution grids. Finally, future work will focus on exploiting the high floating point throughput of GPUs for setting up the IIM stencils, to accelerate the treatment of high order moving boundary problems \cite{Gabbard2024}.

\section*{Acknowledgements}
We wish to acknowledge financial support from an Early Career Award from the Department of Energy, Program Manager Dr.~Steven~Lee, award number DE-SC0020998. 

This research used resources of the National Energy Research Scientific Computing Center (NERSC), a U.S.\ Department of Energy Office of Science User Facility located at Lawrence Berkeley National Laboratory, operated under Contract No. DE-AC02-05CH11231 using NERSC award ASCR-ERCAP0023392.

We thank Gilles Poncelet, Pierre Balty, Thomas Gillis, Dr~Matthieu Duponcheel, and Dr~Philippe Chatelain from UC Louvain for insightful discussions on the Poisson solver and its implementation in MURPHY.

\bibliographystyle{elsarticle-num-names} 
\bibliography{refs.bib}

\begin{thebibliography}{50}
\expandafter\ifx\csname natexlab\endcsname\relax\def\natexlab#1{#1}\fi
\providecommand{\url}[1]{\texttt{#1}}
\providecommand{\href}[2]{#2}
\providecommand{\path}[1]{#1}
\providecommand{\DOIprefix}{doi:}
\providecommand{\ArXivprefix}{arXiv:}
\providecommand{\URLprefix}{URL: }
\providecommand{\Pubmedprefix}{pmid:}
\providecommand{\doi}[1]{\href{http://dx.doi.org/#1}{\path{#1}}}
\providecommand{\Pubmed}[1]{\href{pmid:#1}{\path{#1}}}
\providecommand{\bibinfo}[2]{#2}
\ifx\xfnm\relax \def\xfnm[#1]{\unskip,\space#1}\fi
\bibitem[{Peskin(1972)}]{Peskin1972}
\bibinfo{author}{C.~S. Peskin},
\newblock \bibinfo{title}{{Flow patterns around heart valves: a numerical
  method}},
\newblock \bibinfo{journal}{{Journal of Computational Physics}}
  \bibinfo{volume}{10} (\bibinfo{year}{1972}) \bibinfo{pages}{252--271}.
  \DOIprefix\doi{10.1016/0021-9991(72)90065-4}.
\bibitem[{Shortley and Weller(1938)}]{Shortley1938}
\bibinfo{author}{G.~H. Shortley}, \bibinfo{author}{R.~Weller},
\newblock \bibinfo{title}{The numerical solution of {Laplace}'s equation},
\newblock \bibinfo{journal}{Journal of Applied Physics} \bibinfo{volume}{9}
  (\bibinfo{year}{1938}) \bibinfo{pages}{334--348}. \URLprefix
  \url{https://doi.org/10.1063/1.1710426}. \DOIprefix\doi{10.1063/1.1710426}.
\bibitem[{Jomaa and Macaskill(2005)}]{Jomaa2005}
\bibinfo{author}{Z.~Jomaa}, \bibinfo{author}{C.~Macaskill},
\newblock \bibinfo{title}{The embedded finite difference method for the poisson
  equation in a domain with an irregular boundary and dirichlet boundary
  conditions},
\newblock \bibinfo{journal}{Journal of Computational Physics}
  \bibinfo{volume}{202} (\bibinfo{year}{2005}) \bibinfo{pages}{488–506}.
  \URLprefix \url{http://dx.doi.org/10.1016/j.jcp.2004.07.011}.
  \DOIprefix\doi{10.1016/j.jcp.2004.07.011}.
\bibitem[{Mayo(1984)}]{Mayo1984}
\bibinfo{author}{A.~Mayo},
\newblock \bibinfo{title}{The fast solution of poisson's and the biharmonic
  equations on irregular regions},
\newblock \bibinfo{journal}{SIAM Journal on Numerical Analysis}
  \bibinfo{volume}{21} (\bibinfo{year}{1984}) \bibinfo{pages}{285--299}.
  \URLprefix \url{https://doi.org/10.1137/0721021}.
  \DOIprefix\doi{10.1137/0721021}.
\bibitem[{LeVeque and Li(1994)}]{Leveque1994}
\bibinfo{author}{R.~J. LeVeque}, \bibinfo{author}{Z.~Li},
\newblock \bibinfo{title}{{The immersed interface method for elliptic equations
  with discontinuous coefficients and singular sources}},
\newblock \bibinfo{journal}{SIAM Journal on Numerical Analysis}
  \bibinfo{volume}{31} (\bibinfo{year}{1994}) \bibinfo{pages}{1019--1044}.
  \DOIprefix\doi{10.1137/0731054}.
\bibitem[{Johansen and Colella(1998)}]{johansen1998cartesian}
\bibinfo{author}{H.~Johansen}, \bibinfo{author}{P.~Colella},
\newblock \bibinfo{title}{A cartesian grid embedded boundary method for
  poisson's equation on irregular domains},
\newblock \bibinfo{journal}{Journal of Computational Physics}
  \bibinfo{volume}{147} (\bibinfo{year}{1998}) \bibinfo{pages}{60--85}.
  \DOIprefix\doi{10.1006/jcph.1998.5965}.
\bibitem[{Crockett et~al.(2011)Crockett, Colella, and Graves}]{Crockett2011}
\bibinfo{author}{R.~Crockett}, \bibinfo{author}{P.~Colella},
  \bibinfo{author}{D.~Graves},
\newblock \bibinfo{title}{A cartesian grid embedded boundary method for solving
  the poisson and heat equations with discontinuous coefficients in three
  dimensions},
\newblock \bibinfo{journal}{Journal of Computational Physics}
  \bibinfo{volume}{230} (\bibinfo{year}{2011}) \bibinfo{pages}{2451–2469}.
  \URLprefix \url{http://dx.doi.org/10.1016/j.jcp.2010.12.017}.
  \DOIprefix\doi{10.1016/j.jcp.2010.12.017}.
\bibitem[{Fedkiw et~al.(1999)Fedkiw, Aslam, Merriman, and Osher}]{Fedkiw1999}
\bibinfo{author}{R.~P. Fedkiw}, \bibinfo{author}{T.~Aslam},
  \bibinfo{author}{B.~Merriman}, \bibinfo{author}{S.~Osher},
\newblock \bibinfo{title}{{A non-oscillatory Eulerian approach to interfaces in
  multimaterial flows (the ghost fluid method)}},
\newblock \bibinfo{journal}{{Journal of Computational Physics}}
  \bibinfo{volume}{152} (\bibinfo{year}{1999}) \bibinfo{pages}{457--492}.
  \DOIprefix\doi{10.1006/jcph.1999.6236}.
\bibitem[{Gibou et~al.(2002)Gibou, Fedkiw, Cheng, and Kang}]{Gibou2002}
\bibinfo{author}{F.~Gibou}, \bibinfo{author}{R.~P. Fedkiw},
  \bibinfo{author}{L.-T. Cheng}, \bibinfo{author}{M.~Kang},
\newblock \bibinfo{title}{A second-order-accurate symmetric discretization of
  the poisson equation on irregular domains},
\newblock \bibinfo{journal}{Journal of Computational Physics}
  \bibinfo{volume}{176} (\bibinfo{year}{2002}) \bibinfo{pages}{205–227}.
  \URLprefix \url{http://dx.doi.org/10.1006/jcph.2001.6977}.
  \DOIprefix\doi{10.1006/jcph.2001.6977}.
\bibitem[{de~Prenter et~al.(2023)de~Prenter, Verhoosel, van Brummelen, Larson,
  and Badia}]{dePrenter2023}
\bibinfo{author}{F.~de~Prenter}, \bibinfo{author}{C.~V. Verhoosel},
  \bibinfo{author}{E.~H. van Brummelen}, \bibinfo{author}{M.~G. Larson},
  \bibinfo{author}{S.~Badia},
\newblock \bibinfo{title}{Stability and conditioning of immersed finite element
  methods: Analysis and remedies},
\newblock \bibinfo{journal}{Archives of Computational Methods in Engineering}
  \bibinfo{volume}{30} (\bibinfo{year}{2023}) \bibinfo{pages}{3617–3656}.
  \URLprefix \url{http://dx.doi.org/10.1007/s11831-023-09913-0}.
  \DOIprefix\doi{10.1007/s11831-023-09913-0}.
\bibitem[{Saye(2017)}]{Saye2017}
\bibinfo{author}{R.~Saye},
\newblock \bibinfo{title}{Implicit mesh discontinuous galerkin methods and
  interfacial gauge methods for high-order accurate interface dynamics, with
  applications to surface tension dynamics, rigid body fluid–structure
  interaction, and free surface flow: Part i},
\newblock \bibinfo{journal}{Journal of Computational Physics}
  \bibinfo{volume}{344} (\bibinfo{year}{2017}) \bibinfo{pages}{647–682}.
  \URLprefix \url{http://dx.doi.org/10.1016/j.jcp.2017.04.076}.
  \DOIprefix\doi{10.1016/j.jcp.2017.04.076}.
\bibitem[{Stein et~al.(2016)Stein, Guy, and Thomases}]{Stein2016}
\bibinfo{author}{D.~B. Stein}, \bibinfo{author}{R.~D. Guy},
  \bibinfo{author}{B.~Thomases},
\newblock \bibinfo{title}{Immersed boundary smooth extension: A high-order
  method for solving pde on arbitrary smooth domains using fourier spectral
  methods},
\newblock \bibinfo{journal}{Journal of Computational Physics}
  \bibinfo{volume}{304} (\bibinfo{year}{2016}) \bibinfo{pages}{252–274}.
  \URLprefix \url{http://dx.doi.org/10.1016/j.jcp.2015.10.023}.
  \DOIprefix\doi{10.1016/j.jcp.2015.10.023}.
\bibitem[{Marques et~al.(2017)Marques, Nave, and Rosales}]{Marques2017}
\bibinfo{author}{A.~N. Marques}, \bibinfo{author}{J.-C. Nave},
  \bibinfo{author}{R.~R. Rosales},
\newblock \bibinfo{title}{High order solution of poisson problems with
  piecewise constant coefficients and interface jumps},
\newblock \bibinfo{journal}{Journal of Computational Physics}
  \bibinfo{volume}{335} (\bibinfo{year}{2017}) \bibinfo{pages}{497--515}.
  \URLprefix \url{https://doi.org/10.1016/j.jcp.2017.01.029}.
  \DOIprefix\doi{10.1016/j.jcp.2017.01.029}.
\bibitem[{Ito et~al.(2005)Ito, Li, and Kyei}]{Ito2005}
\bibinfo{author}{K.~Ito}, \bibinfo{author}{Z.~Li}, \bibinfo{author}{Y.~Kyei},
\newblock \bibinfo{title}{Higher-order, cartesian grid based finite difference
  schemes for elliptic equations on irregular domains},
\newblock \bibinfo{journal}{SIAM Journal on Scientific Computing}
  \bibinfo{volume}{27} (\bibinfo{year}{2005}) \bibinfo{pages}{346--367}.
  \URLprefix \url{https://doi.org/10.1137/03060120X}.
  \DOIprefix\doi{10.1137/03060120x}.
\bibitem[{Linnick and Fasel(2005)}]{Linnick2005}
\bibinfo{author}{M.~N. Linnick}, \bibinfo{author}{H.~F. Fasel},
\newblock \bibinfo{title}{{A high-order immersed interface method for
  simulating unsteady incompressible flows on irregular domains}},
\newblock \bibinfo{journal}{{Journal of Computational Physics}}
  \bibinfo{volume}{204} (\bibinfo{year}{2005}) \bibinfo{pages}{157--192}.
  \DOIprefix\doi{10.1016/j.jcp.2004.09.017}.
\bibitem[{Gibou and Fedkiw(2005)}]{Gibou2005}
\bibinfo{author}{F.~Gibou}, \bibinfo{author}{R.~Fedkiw},
\newblock \bibinfo{title}{A fourth order accurate discretization for the
  {Laplace} and heat equations on arbitrary domains, with applications to the
  {Stefan} problem},
\newblock \bibinfo{journal}{{Journal of Computational Physics}}
  \bibinfo{volume}{202} (\bibinfo{year}{2005}) \bibinfo{pages}{577--601}.
  \DOIprefix\doi{10.1016/j.jcp.2004.07.018}.
\bibitem[{Zhou et~al.(2006)Zhou, Zhao, Feig, and Wei}]{Zhou2006}
\bibinfo{author}{Y.~Zhou}, \bibinfo{author}{S.~Zhao},
  \bibinfo{author}{M.~Feig}, \bibinfo{author}{G.~Wei},
\newblock \bibinfo{title}{High order matched interface and boundary method for
  elliptic equations with discontinuous coefficients and singular sources},
\newblock \bibinfo{journal}{Journal of Computational Physics}
  \bibinfo{volume}{213} (\bibinfo{year}{2006}) \bibinfo{pages}{1--30}.
  \URLprefix \url{https://doi.org/10.1016/j.jcp.2005.07.022}.
  \DOIprefix\doi{10.1016/j.jcp.2005.07.022}.
\bibitem[{Zhou and Wei(2006)}]{Zhou2006a}
\bibinfo{author}{Y.~Zhou}, \bibinfo{author}{G.-W. Wei},
\newblock \bibinfo{title}{On the fictitious-domain and interpolation
  formulations of the matched interface and boundary (mib) method},
\newblock \bibinfo{journal}{Journal of Computational Physics}
  \bibinfo{volume}{219} (\bibinfo{year}{2006}) \bibinfo{pages}{228--246}.
\bibitem[{Marques et~al.(2011)Marques, Nave, and Rosales}]{Marques2011}
\bibinfo{author}{A.~N. Marques}, \bibinfo{author}{J.-C. Nave},
  \bibinfo{author}{R.~R. Rosales},
\newblock \bibinfo{title}{A correction function method for poisson problems
  with interface jump conditions},
\newblock \bibinfo{journal}{Journal of Computational Physics}
  \bibinfo{volume}{230} (\bibinfo{year}{2011}) \bibinfo{pages}{7567--7597}.
\bibitem[{Zhu et~al.(2016)Zhu, Luo, and Li}]{Zhu2016}
\bibinfo{author}{C.~Zhu}, \bibinfo{author}{H.~Luo}, \bibinfo{author}{G.~Li},
\newblock \bibinfo{title}{High-order immersed-boundary method for
  incompressible flows},
\newblock \bibinfo{journal}{AIAA Journal} \bibinfo{volume}{54}
  (\bibinfo{year}{2016}) \bibinfo{pages}{2734--2741}.
  \DOIprefix\doi{10.2514/1.J054628}.
\bibitem[{Hosseinverdi and Fasel(2018)}]{Hosseinverdi2018}
\bibinfo{author}{S.~Hosseinverdi}, \bibinfo{author}{H.~F. Fasel},
\newblock \bibinfo{title}{An efficient, high-order method for solving {Poisson}
  equation for immersed boundaries: {Combination} of compact difference and
  multiscale multigrid methods},
\newblock \bibinfo{journal}{{Journal of Computational Physics}}
  \bibinfo{volume}{374} (\bibinfo{year}{2018}) \bibinfo{pages}{912--940}.
  \DOIprefix\doi{10.1016/j.jcp.2018.08.006}.
\bibitem[{Hosseinverdi and Fasel(2020)}]{Hosseinverdi2020}
\bibinfo{author}{S.~Hosseinverdi}, \bibinfo{author}{H.~F. Fasel},
\newblock \bibinfo{title}{A fourth-order accurate compact difference scheme for
  solving the three-dimensional {Poisson} equation with arbitrary boundaries},
\newblock in: \bibinfo{booktitle}{AIAA Scitech 2020 Forum},
  \bibinfo{year}{2020}, p. \bibinfo{pages}{0805}.
  \DOIprefix\doi{10.2514/6.2020-0805}.
\bibitem[{Devendran et~al.(2017)Devendran, Graves, Johansen, and
  Ligocki}]{devendran2017fourth}
\bibinfo{author}{D.~Devendran}, \bibinfo{author}{D.~Graves},
  \bibinfo{author}{H.~Johansen}, \bibinfo{author}{T.~Ligocki},
\newblock \bibinfo{title}{A fourth-order {Cartesian} grid embedded boundary
  method for {Poisson’s} equation},
\newblock \bibinfo{journal}{Communications in Applied Mathematics and
  Computational Science} \bibinfo{volume}{12} (\bibinfo{year}{2017})
  \bibinfo{pages}{51--79}. \DOIprefix\doi{10.2140/camcos.2017.12.51}.
\bibitem[{Thacher et~al.(2023)Thacher, Johansen, and Martin}]{Thacher2023AHigh}
\bibinfo{author}{W.~Thacher}, \bibinfo{author}{H.~Johansen},
  \bibinfo{author}{D.~Martin},
\newblock \bibinfo{title}{A high order cartesian grid, finite volume method for
  elliptic interface problems},
\newblock \bibinfo{journal}{Journal of Computational Physics}
  \bibinfo{volume}{491} (\bibinfo{year}{2023}) \bibinfo{pages}{112351}.
  \URLprefix \url{https://doi.org/10.1016/j.jcp.2023.112351}.
  \DOIprefix\doi{10.1016/j.jcp.2023.112351}.
\bibitem[{Overton-Katz et~al.(2023)Overton-Katz, Gao, Guzik, Antepara, Graves,
  and Johansen}]{OvertonKatz2022a}
\bibinfo{author}{N.~Overton-Katz}, \bibinfo{author}{X.~Gao},
  \bibinfo{author}{S.~Guzik}, \bibinfo{author}{O.~Antepara},
  \bibinfo{author}{D.~T. Graves}, \bibinfo{author}{H.~Johansen},
\newblock \bibinfo{title}{A fourth-order embedded boundary finite volume method
  for the unsteady stokes equations with complex geometries},
\newblock \bibinfo{journal}{SIAM Journal on Scientific Computing}
  \bibinfo{volume}{45} (\bibinfo{year}{2023}) \bibinfo{pages}{A2409--A2430}.
  \DOIprefix\doi{10.1137/22M1532019}.
\bibitem[{Thacher et~al.(2024)Thacher, Johansen, and Martin}]{thacher2024high}
\bibinfo{author}{W.~Thacher}, \bibinfo{author}{H.~Johansen},
  \bibinfo{author}{D.~Martin},
\newblock \bibinfo{title}{A high order cut-cell method for solving the
  shallow-shelf equations},
\newblock \bibinfo{journal}{Journal of Computational Science}
  \bibinfo{volume}{80} (\bibinfo{year}{2024}) \bibinfo{pages}{102319}.
\bibitem[{Devendran et~al.(2014)Devendran, Graves, and
  Johansen}]{Devendran2014}
\bibinfo{author}{D.~Devendran}, \bibinfo{author}{D.~T. Graves},
  \bibinfo{author}{H.~Johansen}, \bibinfo{title}{A Hybrid Multigrid Algorithm
  for Poisson's Equation using an Adaptive, Fourth Order Treatment of Cut
  Cells}, \bibinfo{type}{Technical Report} \bibinfo{number}{LBNL-1004329},
  Lawrence Berkeley National Laboratory, \bibinfo{year}{2014}. \URLprefix
  \url{https://crd.lbl.gov/assets/pubs_presos/multigrid.pdf}.
\bibitem[{Adams and Li(2002)}]{Adams:2002}
\bibinfo{author}{L.~Adams}, \bibinfo{author}{Z.~Li},
\newblock \bibinfo{title}{The immersed interface/multigrid methods for
  interface problems},
\newblock \bibinfo{journal}{SIAM Journal on Scientific Computing}
  \bibinfo{volume}{24} (\bibinfo{year}{2002}) \bibinfo{pages}{463–479}.
  \URLprefix \url{http://dx.doi.org/10.1137/S1064827501389849}.
  \DOIprefix\doi{10.1137/s1064827501389849}.
\bibitem[{Adams and Chartier(2004)}]{Adams:2004}
\bibinfo{author}{L.~Adams}, \bibinfo{author}{T.~P. Chartier},
\newblock \bibinfo{title}{New geometric immersed interface multigrid solvers},
\newblock \bibinfo{journal}{SIAM Journal on Scientific Computing}
  \bibinfo{volume}{25} (\bibinfo{year}{2004}) \bibinfo{pages}{1516–1533}.
  \URLprefix \url{http://dx.doi.org/10.1137/S1064827503421707}.
  \DOIprefix\doi{10.1137/s1064827503421707}.
\bibitem[{Chen and Strain(2008)}]{Chen:2008}
\bibinfo{author}{T.~Chen}, \bibinfo{author}{J.~Strain},
\newblock \bibinfo{title}{Piecewise-polynomial discretization and
  krylov-accelerated multigrid for elliptic interface problems},
\newblock \bibinfo{journal}{Journal of Computational Physics}
  \bibinfo{volume}{227} (\bibinfo{year}{2008}) \bibinfo{pages}{7503–7542}.
  \URLprefix \url{http://dx.doi.org/10.1016/j.jcp.2008.04.027}.
  \DOIprefix\doi{10.1016/j.jcp.2008.04.027}.
\bibitem[{Gillis et~al.(2018)Gillis, Winckelmans, and Chatelain}]{Gillis2018}
\bibinfo{author}{T.~Gillis}, \bibinfo{author}{G.~Winckelmans},
  \bibinfo{author}{P.~Chatelain},
\newblock \bibinfo{title}{Fast immersed interface poisson solver for {3D}
  unbounded problems around arbitrary geometries},
\newblock \bibinfo{journal}{Journal of Computational Physics}
  \bibinfo{volume}{354} (\bibinfo{year}{2018}) \bibinfo{pages}{403--416}.
  \URLprefix \url{https://doi.org/10.1016/j.jcp.2017.10.042}.
  \DOIprefix\doi{10.1016/j.jcp.2017.10.042}.
\bibitem[{Gabbard and van Rees(2023)}]{Gabbard2023}
\bibinfo{author}{J.~Gabbard}, \bibinfo{author}{W.~M. van Rees},
\newblock \bibinfo{title}{A high-order {3D} immersed interface finite
  difference method for the advection-diffusion equation},
\newblock in: \bibinfo{booktitle}{AIAA SCITECH 2023 Forum},
  \bibinfo{year}{2023}, p. \bibinfo{pages}{2480}.
  \DOIprefix\doi{10.2514/6.2023-2480}.
\bibitem[{Gabbard and van Rees(2024)}]{Gabbard2024}
\bibinfo{author}{J.~Gabbard}, \bibinfo{author}{W.~M. van Rees},
\newblock \bibinfo{title}{A high-order finite difference method for moving
  immersed domain boundaries and material interfaces},
\newblock \bibinfo{journal}{Journal of Computational Physics}
  \bibinfo{volume}{507} (\bibinfo{year}{2024}) \bibinfo{pages}{112979}.
  \URLprefix \url{https://doi.org/10.1016/j.jcp.2024.112979}.
  \DOIprefix\doi{10.1016/j.jcp.2024.112979}.
\bibitem[{Li and Ito(2006)}]{Li2006}
\bibinfo{author}{Z.~Li}, \bibinfo{author}{K.~Ito}, \bibinfo{title}{The immersed
  interface method: numerical solutions of PDEs involving interfaces and
  irregular domains}, \bibinfo{publisher}{SIAM}, \bibinfo{year}{2006}.
\bibitem[{Wiegmann and Bube(2000)}]{Wiegmann2000}
\bibinfo{author}{A.~Wiegmann}, \bibinfo{author}{K.~P. Bube},
\newblock \bibinfo{title}{{The explicit-jump immersed interface method: finite
  difference methods for PDEs with piecewise smooth solutions}},
\newblock \bibinfo{journal}{SIAM Journal on Numerical Analysis}
  \bibinfo{volume}{37} (\bibinfo{year}{2000}) \bibinfo{pages}{827 -- 862}.
\bibitem[{Hosseinverdi and Fasel(2017)}]{Hosseinverdi2017}
\bibinfo{author}{S.~Hosseinverdi}, \bibinfo{author}{H.~F. Fasel},
\newblock \bibinfo{title}{Very high-order accurate sharp immersed interface
  method: application to direct numerical simulations of incompressible flows},
\newblock in: \bibinfo{booktitle}{23rd AIAA Computational Fluid Dynamics
  Conference}, \bibinfo{year}{2017}, p. \bibinfo{pages}{3624}.
  \DOIprefix\doi{10.2514/6.2017-3624}.
\bibitem[{Gabbard et~al.(2022)Gabbard, Gillis, Chatelain, and van
  Rees}]{Gabbard2022}
\bibinfo{author}{J.~Gabbard}, \bibinfo{author}{T.~Gillis},
  \bibinfo{author}{P.~Chatelain}, \bibinfo{author}{W.~M. van Rees},
\newblock \bibinfo{title}{An immersed interface method for the {2D}
  vorticity-velocity {Navier--Stokes} equations with multiple bodies},
\newblock \bibinfo{journal}{{Journal of Computational Physics}}
  \bibinfo{volume}{464} (\bibinfo{year}{2022}) \bibinfo{pages}{111339}.
  \DOIprefix\doi{10.1016/j.jcp.2022.111339}.
\bibitem[{Li(1997)}]{Li1997}
\bibinfo{author}{Z.~Li},
\newblock \bibinfo{title}{Immersed interface methods for moving interface
  problems},
\newblock \bibinfo{journal}{Numerical algorithms} \bibinfo{volume}{14}
  (\bibinfo{year}{1997}) \bibinfo{pages}{269--293}.
  \DOIprefix\doi{10.1023/A:1019173215885}.
\bibitem[{Trottenberg et~al.(2000)Trottenberg, Oosterlee, and
  Schuller}]{Trottenberg2000}
\bibinfo{author}{U.~Trottenberg}, \bibinfo{author}{C.~W. Oosterlee},
  \bibinfo{author}{A.~Schuller}, \bibinfo{title}{Multigrid},
  \bibinfo{publisher}{Elsevier}, \bibinfo{year}{2000}.
\bibitem[{Gallinato and Poignard(2015)}]{Gallinato2015}
\bibinfo{author}{O.~Gallinato}, \bibinfo{author}{C.~Poignard},
  \bibinfo{title}{Superconvergent Cartesian Methods for Poisson type Equations
  in {2D}–domains}, \bibinfo{type}{Technical Report} \bibinfo{number}{8809},
  INRIA Bordeaux, \bibinfo{year}{2015}.
\bibitem[{Yoon and Min(2015)}]{Yoon2015}
\bibinfo{author}{G.~Yoon}, \bibinfo{author}{C.~Min},
\newblock \bibinfo{title}{Analyses on the finite difference method by gibou et
  al. for poisson equation},
\newblock \bibinfo{journal}{Journal of Computational Physics}
  \bibinfo{volume}{280} (\bibinfo{year}{2015}) \bibinfo{pages}{184–194}.
  \URLprefix \url{http://dx.doi.org/10.1016/j.jcp.2014.09.009}.
  \DOIprefix\doi{10.1016/j.jcp.2014.09.009}.
\bibitem[{Gillis and van Rees(2022)}]{Gillis2022}
\bibinfo{author}{T.~Gillis}, \bibinfo{author}{W.~M. van Rees},
\newblock \bibinfo{title}{{MURPHY}---a scalable multiresolution framework for
  scientific computing on {3D} block-structured collocated grids},
\newblock \bibinfo{journal}{SIAM Journal on Scientific Computing}
  \bibinfo{volume}{44} (\bibinfo{year}{2022}) \bibinfo{pages}{C367--C398}.
  \DOIprefix\doi{10.1137/21M141676X}.
\bibitem[{Matsunaga and Yamamoto(2000)}]{Matsunaga2000}
\bibinfo{author}{N.~Matsunaga}, \bibinfo{author}{T.~Yamamoto},
\newblock \bibinfo{title}{Superconvergence of the shortley--weller
  approximation for dirichlet problems},
\newblock \bibinfo{journal}{Journal of computational and applied mathematics}
  \bibinfo{volume}{116} (\bibinfo{year}{2000}) \bibinfo{pages}{263--273}.
\bibitem[{Burstedde et~al.(2011)Burstedde, Wilcox, and
  Ghattas}]{Burstedde:2011}
\bibinfo{author}{C.~Burstedde}, \bibinfo{author}{L.~Wilcox},
  \bibinfo{author}{O.~Ghattas},
\newblock \bibinfo{title}{p4est: Scalable algorithms for parallel adaptive mesh
  refinement on forests of octrees},
\newblock \bibinfo{journal}{SIAM Journal on Scientific Computing}
  \bibinfo{volume}{33} (\bibinfo{year}{2011}) \bibinfo{pages}{1103--1133}.
  \URLprefix \url{https://doi.org/10.1137/100791634}.
  \DOIprefix\doi{10.1137/100791634}.
\bibitem[{Balay et~al.(2024)Balay, Abhyankar, Adams, Benson, Brown, Brune,
  Buschelman, Constantinescu, Dalcin, Dener, Eijkhout, Faibussowitsch, Gropp,
  Hapla, Isaac, Jolivet, Karpeev, Kaushik, Knepley, Kong, Kruger, May, McInnes,
  Mills, Mitchell, Munson, Roman, Rupp, Sanan, Sarich, Smith, Zampini, Zhang,
  Zhang, and Zhang}]{petsc-web-page}
\bibinfo{author}{S.~Balay}, \bibinfo{author}{S.~Abhyankar},
  \bibinfo{author}{M.~F. Adams}, \bibinfo{author}{S.~Benson},
  \bibinfo{author}{J.~Brown}, \bibinfo{author}{P.~Brune},
  \bibinfo{author}{K.~Buschelman}, \bibinfo{author}{E.~M. Constantinescu},
  \bibinfo{author}{L.~Dalcin}, \bibinfo{author}{A.~Dener},
  \bibinfo{author}{V.~Eijkhout}, \bibinfo{author}{J.~Faibussowitsch},
  \bibinfo{author}{W.~D. Gropp}, \bibinfo{author}{V.~Hapla},
  \bibinfo{author}{T.~Isaac}, \bibinfo{author}{P.~Jolivet},
  \bibinfo{author}{D.~Karpeev}, \bibinfo{author}{D.~Kaushik},
  \bibinfo{author}{M.~G. Knepley}, \bibinfo{author}{F.~Kong},
  \bibinfo{author}{S.~Kruger}, \bibinfo{author}{D.~A. May},
  \bibinfo{author}{L.~C. McInnes}, \bibinfo{author}{R.~T. Mills},
  \bibinfo{author}{L.~Mitchell}, \bibinfo{author}{T.~Munson},
  \bibinfo{author}{J.~E. Roman}, \bibinfo{author}{K.~Rupp},
  \bibinfo{author}{P.~Sanan}, \bibinfo{author}{J.~Sarich},
  \bibinfo{author}{B.~F. Smith}, \bibinfo{author}{S.~Zampini},
  \bibinfo{author}{H.~Zhang}, \bibinfo{author}{H.~Zhang},
  \bibinfo{author}{J.~Zhang}, \bibinfo{title}{{PETS}c {W}eb page},
  \bibinfo{howpublished}{\url{https://petsc.org/}}, \bibinfo{year}{2024}.
  \URLprefix \url{https://petsc.org/}.
\bibitem[{Poncelet et~al.(2025)Poncelet, Lambrechts, Gillis, and
  Chatelain}]{Poncelet2025}
\bibinfo{author}{G.~Poncelet}, \bibinfo{author}{J.~Lambrechts},
  \bibinfo{author}{T.~Gillis}, \bibinfo{author}{P.~Chatelain},
  \bibinfo{title}{A high order adaptive multiresolution {Poisson} solver for
  problems with unbounded directions}, \bibinfo{year}{2025}. \bibinfo{note}{In
  preparation}.
\bibitem[{Balay et~al.(2024)Balay, Abhyankar, Adams, Benson, Brown, Brune,
  Buschelman, Constantinescu, Dalcin, Dener, Eijkhout, Faibussowitsch, Gropp,
  Hapla, Isaac, Jolivet, Karpeev, Kaushik, Knepley, Kong, Kruger, May, McInnes,
  Mills, Mitchell, Munson, Roman, Rupp, Sanan, Sarich, Smith, Suh, Zampini,
  Zhang, Zhang, and Zhang}]{petsc-user-ref}
\bibinfo{author}{S.~Balay}, \bibinfo{author}{S.~Abhyankar},
  \bibinfo{author}{M.~F. Adams}, \bibinfo{author}{S.~Benson},
  \bibinfo{author}{J.~Brown}, \bibinfo{author}{P.~Brune},
  \bibinfo{author}{K.~Buschelman}, \bibinfo{author}{E.~Constantinescu},
  \bibinfo{author}{L.~Dalcin}, \bibinfo{author}{A.~Dener},
  \bibinfo{author}{V.~Eijkhout}, \bibinfo{author}{J.~Faibussowitsch},
  \bibinfo{author}{W.~D. Gropp}, \bibinfo{author}{V.~Hapla},
  \bibinfo{author}{T.~Isaac}, \bibinfo{author}{P.~Jolivet},
  \bibinfo{author}{D.~Karpeev}, \bibinfo{author}{D.~Kaushik},
  \bibinfo{author}{M.~G. Knepley}, \bibinfo{author}{F.~Kong},
  \bibinfo{author}{S.~Kruger}, \bibinfo{author}{D.~A. May},
  \bibinfo{author}{L.~C. McInnes}, \bibinfo{author}{R.~T. Mills},
  \bibinfo{author}{L.~Mitchell}, \bibinfo{author}{T.~Munson},
  \bibinfo{author}{J.~E. Roman}, \bibinfo{author}{K.~Rupp},
  \bibinfo{author}{P.~Sanan}, \bibinfo{author}{J.~Sarich},
  \bibinfo{author}{B.~F. Smith}, \bibinfo{author}{H.~Suh},
  \bibinfo{author}{S.~Zampini}, \bibinfo{author}{H.~Zhang},
  \bibinfo{author}{H.~Zhang}, \bibinfo{author}{J.~Zhang},
  \bibinfo{title}{{PETSc/TAO} Users Manual}, \bibinfo{type}{Technical Report}
  \bibinfo{number}{ANL-21/39 - Revision 3.22}, Argonne National Laboratory,
  \bibinfo{year}{2024}. \DOIprefix\doi{10.2172/2205494}.
\bibitem[{Saad(2003)}]{Saad:2003}
\bibinfo{author}{Y.~Saad}, \bibinfo{title}{Iterative Methods for Sparse Linear
  Systems}, \bibinfo{edition}{2nd} ed., \bibinfo{publisher}{Society for
  Industrial and Applied Mathematics}, \bibinfo{year}{2003}. \URLprefix
  \url{https://epubs.siam.org/doi/abs/10.1137/1.9780898718003}.
  \DOIprefix\doi{10.1137/1.9780898718003}.
\bibitem[{Saad(1993)}]{Saad:1993}
\bibinfo{author}{Y.~Saad},
\newblock \bibinfo{title}{A {{Flexible Inner-Outer Preconditioned GMRES
  Algorithm}}},
\newblock \bibinfo{journal}{SIAM Journal on Scientific Computing}
  \bibinfo{volume}{14} (\bibinfo{year}{1993}) \bibinfo{pages}{461--469}.
  \DOIprefix\doi{10.1137/0914028}.
\bibitem[{Min and Gibou(2006)}]{Min2006}
\bibinfo{author}{C.~Min}, \bibinfo{author}{F.~Gibou},
\newblock \bibinfo{title}{A second order accurate projection method for the
  incompressible navier–stokes equations on non-graded adaptive grids},
\newblock \bibinfo{journal}{Journal of Computational Physics}
  \bibinfo{volume}{219} (\bibinfo{year}{2006}) \bibinfo{pages}{912–929}.
  \URLprefix \url{http://dx.doi.org/10.1016/j.jcp.2006.07.019}.
  \DOIprefix\doi{10.1016/j.jcp.2006.07.019}.

\end{thebibliography}

\end{document}